\documentclass[a4paper,12pt]{amsart}
\usepackage[centertags]{amsmath}
\usepackage{amsfonts}
\usepackage{amssymb}
\usepackage{amsthm}
\usepackage{newlfont}
\usepackage{color}
\usepackage{graphicx}
\usepackage{tikz}

\def\ee{\end{equation}}
 \theoremstyle{definition}

 \numberwithin{equation}{section}
 \newcommand{\A}{\mathcal{A}}

\def \R {\mathbb{R}}

\def \C {\mathbb{C}}

\def \A {\mathcal{A}}

\newcommand{\delim}[3]{\left#1 #3 \right#2}

\newcommand{\norma}[1]{\delim{\|}{\|}{#1}}

\def \ds {\displaystyle}

\newtheorem{theorem}{Theorem}[section]

\newtheorem{lemma}[theorem]{Lemma}

\newtheorem{proposition}[theorem]{Proposition}

\newtheorem{definition}[theorem]{Definition}

\newtheorem{remark}[theorem]{Remark}

\title[Dynamics in Dumbbell domains]{Pullback and uniform attractors for nonautonomous reaction-diffusion equation in Dumbbell domains}
\vspace{8pt}

\author[M. Belluzi]{Maykel Belluzi$^\ddag$}\thanks{$^\ddag$Research partially
	supported by FAPESP \# 2017/09406-0 and \# 2017/17502-0, Brazil}
\address[M. Belluzi]{Universidade Federal de S\~{a}o
	Carlos, Departamento de Matem\'atica, 13565-905 S\~{a}o
	Carlos SP, Brazil.}
\email{maykel@dm.ufscar.br}

\author[T. Caraballo]{Tom\'as Caraballo$^\dag$}\thanks{$^\dag$Research partially supported by Ministerio de Ciencia Innovaci\'on y Universidades (Spain), FEDER (European Community) under grant PGC2018- 096540-B-I00, and by Junta de Andaluc\'ia (Consejer\'ia de Econom\'ia y Conocimiento) under project US-1254251.}
\address[T. Caraballo]{Departamento de Ecuaciones Diferenciales y An\'alisis Num\'erico, Universidad de Sevilla, Apdo. de Correos 1160, 41080-Sevilla, Spain.}
\email{caraballo@us.es}

\author[M. J. D. Nascimento]{Marcelo J. D. Nascimento$^\star$}\thanks{$^\star$Research partially
	supported by FAPESP \# 2017/06582-2, Brazil}
\address[M. J. D. Nascimento]{Universidade Federal de S\~{a}o
	Carlos, Departamento de Matem\'atica, 13565-905 S\~{a}o
	Carlos SP, Brazil.}
\email{marcelo@dm.ufscar.br}

\author[K. Schiabel]{Karina Schiabel}
\address[K. Schiabel]{Universidade Federal de S\~{a}o
	Carlos, Departamento de Matem\'atica, 13565-905 S\~{a}o
	Carlos SP, Brazil.}
\email{schiabel@dm.ufscar.br}

\date{\today}

\begin{document}

%%%%%%%%%%%%%%% ABSTRACT %%%%%%%%%%%%%%%

\begin{abstract}
This work is devoted to the study of the asymptotic behavior of  nonautonomous reaction-diffusion equations in Dumbbell domains $\Omega_{\varepsilon} \subset \R^{N}$. Each $\Omega_{\varepsilon}$ is the union of a fixed open set $\Omega$ and a channel $R_{\varepsilon}$ that collapses to a line segment $R_0$ as $\varepsilon \rightarrow 0^{+}$. 
%, which consist of an open set $\Omega \subset R^{N}$ and a channel $\R_{\varepsilon}$ that collapses to a line segment $R_0$
 We first establish the global existence of solution for each problem by using two properties of the parabolic equation considered, which are the positivity of the solutions and comparison results for them. We prove the existence of pullback and uniform attractors and we obtain uniform bounds (in $\varepsilon$) for them.


\medskip

 
\end{abstract}

 \maketitle

 %%%%%%%%%%%INTRODUCTION %%%%%%%%%%%%%%%

\section{Introduction}

In this paper we address the problem of the asymptotic behavior of a nonautonomous nonlinear reaction-diffusion equation in domains of Dumbbell type. Those Dumbbell domains will be denoted by $\Omega_{\varepsilon}$, $\varepsilon \in [0,1]$.

Each $\Omega_{\varepsilon} \subset \R^{N}$, $N \geq 2$, is a smooth domain, consisting of two disconnected sets $\Omega$ joined by a thin channel $R_{\varepsilon}$, which degenerates to a line segment $R_0$ as $\varepsilon \rightarrow 0^{+}$, see Fig. \ref{Fig_1} below. For instance, the channels $R_{\varepsilon}$ are obtained from the channel $R_1$ as follows $R_{\varepsilon} = \{   (s,\varepsilon x') : (s,x') \in R_1 \}$, and $R_1=\{ (s,x') : 0\leq s \leq 1 \mbox{ and } x'\in \Gamma_s^{1}     \}$ where $\Gamma_s^1 \subset \R^{N-1}$ is $\mathcal{C}^1\mbox{-diffeomorphic}$ to the unitary ball of $\R^{N-1}.$

\definecolor{NColor}{rgb}{0.5,0.5,0.5} %Para mudar a cor do desenho, basta mudas os par\^ametros na \'ultima chave
\begin{figure}[!ht]\label{Fig_1}
\[
\begin{minipage}{6.5cm}
\begin{tikzpicture}
\draw[line width=0.8pt,color=NColor] (2,2)  to [out=-10,in=190]  (4,2);
\draw[line width=0.8pt,color=NColor] (2,1.7)  to [out=5,in=175]  (4,1.7);
\draw[line width=0.8pt,color=NColor] (4,1.7)  to [out=260,in=180]  (4.5,0);
\draw[line width=0.8pt,color=NColor] (4.5,0)  to [out=0,in=180]  (5,0.2);
\draw[line width=0.8pt,color=NColor] (5,0.2)  to [out=0,in=180]  (5.5,0);
\draw[line width=0.8pt,color=NColor] (5.5,0)  to [out=0,in=270]  (5.5,2);
\draw[line width=0.8pt,color=NColor] (5.5,2)  to [out=90,in=-20]  (5.5,4);
\draw[line width=0.8pt,color=NColor] (5.5,4)  to [out=160,in=70]  (3.5,3.5);
\draw[line width=0.8pt,color=NColor] (3.5,3.5)  to [out=250,in=100]  (4,2);
\draw[line width=0.8pt,color=NColor] (2,2)  to [out=60,in=0]  (1.2,3.5);
\draw[line width=0.8pt,color=NColor] (1.2,3.5)  to [out=180,in=90]  (0,2);
\draw[line width=0.8pt,color=NColor] (0,2)  to [out=270,in=200]  (2,0);
\draw[line width=0.8pt,color=NColor] (2,0)  to [out=20,in=290]  (2,1.7);
\fill[color=NColor!20] (2,2)  to [out=60,in=0]  (1.2,3.5) to [out=180,in=90]  (0,2) to [out=270,in=200]  (2,0) to [out=20,in=290]  (2,1.7);
\fill[color=NColor!20] (2,2)  to [out=-10,in=190]  (4,2) -- (4,1.7)   to [out=175,in=5]  (2,1.7);
\fill[color=NColor!20] (4,1.7)  to [out=260,in=180]  (4.5,0)  to [out=0,in=180]  (5,0.2)  to [out=0,in=180]  (5.5,0) to [out=0,in=270]  (5.5,2) to [out=90,in=-20]  (5.5,4) to [out=160,in=70]  (3.5,3.5) to [out=250,in=100]  (4,2);
\node at (1.5,0.5) {\footnotesize $\Omega$};
\node at (5,0.5) {\footnotesize $\Omega$};
\node at (3,1.5) {\footnotesize \rm{R}$_{\varepsilon}$};
\node at (2.8,4) {\footnotesize $\Omega_{\varepsilon}$};
\end{tikzpicture}
\end{minipage}
\begin{minipage}{2cm}
\[
\stackrel{\varepsilon \to 0^+}{\tikz \draw[-latex,line width=0.7pt] (0,0) -- (1,0);}
\]
\end{minipage}
\begin{minipage}{6.5cm}
\begin{tikzpicture}
\draw[line width=0.5pt,color=NColor] (1.95,1.85)  --  (4,1.85);
\draw[line width=0.8pt,color=NColor] (4,1.7)  to [out=260,in=180]  (4.5,0);
\draw[line width=0.8pt,color=NColor] (4.5,0)  to [out=0,in=180]  (5,0.2);
\draw[line width=0.8pt,color=NColor] (5,0.2)  to [out=0,in=180]  (5.5,0);
\draw[line width=0.8pt,color=NColor] (5.5,0)  to [out=0,in=270]  (5.5,2);
\draw[line width=0.8pt,color=NColor] (5.5,2)  to [out=90,in=-20]  (5.5,4);
\draw[line width=0.8pt,color=NColor] (5.5,4)  to [out=160,in=70]  (3.5,3.5);
\draw[line width=0.8pt,color=NColor] (3.5,3.5)  to [out=250,in=100]  (4,2) to [out=280,in=80] (4,1.7);
\draw[line width=0.8pt,color=NColor] (2,1.7) to [out=110,in=240] (2,2)  to [out=60,in=0]  (1.2,3.5);
\draw[line width=0.8pt,color=NColor] (1.2,3.5)  to [out=180,in=90]  (0,2);
\draw[line width=0.8pt,color=NColor] (0,2)  to [out=270,in=200]  (2,0);
\draw[line width=0.8pt,color=NColor] (2,0)  to [out=20,in=290]  (2,1.7);
\fill[color=NColor!20] (2,2)  to [out=60,in=0]  (1.2,3.5) to [out=180,in=90]  (0,2) to [out=270,in=200]  (2,0) to [out=20,in=290]  (2,1.7) to [out=110,in=240] (2,2);
\fill[color=NColor!20] (4,1.7)  to [out=260,in=180]  (4.5,0)  to [out=0,in=180]  (5,0.2)  to [out=0,in=180]  (5.5,0) to [out=0,in=270]  (5.5,2) to [out=90,in=-20]  (5.5,4) to [out=160,in=70]  (3.5,3.5) to [out=250,in=100]  (4,2) to [out=280,in=80] (4,1.7);
\node at (1.5,0.5) {\footnotesize $\Omega$};
\node at (5,0.5) {\footnotesize $\Omega$};
\node at (3,1.6) {\footnotesize \rm{R}$_{0}$};
\end{tikzpicture}
\end{minipage}
\]
\caption{\small Dumbbell Domain}
\end{figure}

For each $\varepsilon \in (0,1]$, we consider the nonautonomous evolution equation
\begin{equation}\label{Pe*}
\begin{cases}
u_t-a(t)\Delta u + u = f(u), \quad t> \tau, x\in \Omega_{\varepsilon}, \\
\frac{\partial u}{\partial n} = 0, \qquad \qquad \qquad  \mbox{ }\qquad \partial \Omega_{\varepsilon},
\end{cases}
\end{equation}where $a : \R \rightarrow [c_0,c_1]$, $0<c_0<c_1$ and $f: \R \rightarrow \R^{}$ are continuously differentiable. Furthermore, $a$ is H\"{o}lder continuous with exponent $\delta \in (0,1] $ and constant $C$, that is, 
\begin{equation}\label{a_holder}
|a(t) - a(s)| \leq C|s-t|^{\delta}, \quad \forall t,s\in \R,
\end{equation}
%,com $U$ aberto do $\R^N$ contendo todos os conjuntos $\Omega_{\varepsilon}$ paramentrizados em $\varepsilon$,
and the nonlinearity $f$ satisfies a growth condition
\begin{equation}\label{G**}
|f'(r)| \leq C (1+|r|^{\rho -1}), \quad \mbox{for some }\rho \geq 1.
\end{equation}

As $\varepsilon \rightarrow 0^{+}$, the limit equation obtained is 
\begin{equation}\label{P0*}
\begin{cases}
w_t-a(t) \Delta w + w = f(w), \quad \quad t> \tau, x\in \Omega, \\
\dfrac{\partial w}{\partial n}=0, \quad \hspace{3.8cm} x \in \partial \Omega, \\
v_t-a(t) \dfrac{1}{g} (gv')' + v = f(v), \hspace{0.4cm} t>\tau, s\in (0,1), \\
v(p_0) = w(p_0) \mbox{ and }v(p_1) = w(p_1), \\
\end{cases}
\end{equation}
where $p_0 = (0,0,...,0) \in \R^N$, $p_1 = (1,0,...,0) \in \R^N$ and $g$ is a function that appears as consequence of the geometry of the initial channel $R_1$. If $w \in W^{2,p}(\Omega)$,
%, we will have $w$ with enough regularity ($w\in \mathcal{C}(\overline{\Omega})$) 
then it is continuously extended to $\overline{\Omega}$ and the boundary conditions for $v$  make sense. Note that the dynamics in $\Omega$ is independent of the dynamics in $R_0$, whereas the dynamics in $R_0$ depends on the dynamics in $\Omega$ through the continuity conditions at $p_0$ and $p_1$. We express this feature by saying that this system is one-sided coupled. For further details about limit equations on thin domains, we recommend \cite{HaleRaugel1992,Oliva,Prizzi}.

\color{black}
Those Dumbbell domains emerged in the literature as a counterpart of convex domains in the following sense: if we consider an autonomous reaction-diffusion equation in a convex domain, the stable equilibria for the equation are constant in the domain (see, for instance, \cite{Casten,KKW}). In order to obtain stable equilibria that are not spatially constant, we can not allow the domain to be convex. This is when the Dumbbell domains appeared as prototype of nonconvex domain. The earlier works dealing with Dumbbell domains considered the case in which the channel $R_{\varepsilon}$ was a straight cylinder. In \cite{JimboIII, JimboI, JimboII}, Jimbo made a detailed study of semilinear elliptic equations on those Dumbbell domains with straight cylinder channel and exhibited stable equilibrium that is not spatially constant. The same author proved in \cite{Jimbo_attractor} a result on the existence of global attractor $\A_{\varepsilon}$ for a parabolic semilinear equation in a Dumbbell domain $\Omega_{\varepsilon}$ with straight channel, including the limiting case %, that is,
% the existence of global attractor $\A_0$ for the semilinear equation in the singular domain 
%the equation 
in $\Omega_0.$

Allowing the channel $R_{\varepsilon}$ to be more generic than the cylinders considered by Jimbo, the authors in \cite{DDDI,DDDII,DDDIII} studied an autonomous version of the equations \eqref{Pe*} and \eqref{P0*}, when $a(t) = a>0$. To treat this equations in Dumbbell domains, they developed an appropriate functional setting. 	For $\varepsilon \in (0,1]$, the phase space considered was $U_p^{\varepsilon}:= L^{p}  (\Omega_{\varepsilon})$, with the norm 
\begin{equation}\label{norma_e}
\norma{u}_{U_p^{\varepsilon}} = \norma{u}_{L^{p}(\Omega)}   + \frac{1}{\varepsilon^{\frac{N-1}{p}}}  \norma{u}_{L^{p}(R_{\varepsilon})} \end{equation}
and, for $\varepsilon =0$, $U_p^{0}:=L^{p}(\Omega) \oplus   L^{p}_g(0,1)$,  that is $(w,v) \in U_p^{0}$ if $w\in L^{p}(\Omega)$, $v\in L^{p}(0,1)$ and is equipped with norm
\begin{equation}\label{norma_0}
\norma{(w,v)}_{U_p^{0}}   = \norma{w}_{L^{p}(\Omega)}   +\left[   \int_0^{1} g(s) |v(s)|^{p}ds  \right]^{\frac{1}{p}}. 
\end{equation} 

	This functional setting allowed them to prove the existence of stable equilibrium that was not constant in the domain $\Omega_{\varepsilon}$. The construction of such equilibrium involved a detailed study of the limiting problem in the domain $\Omega_0 = \Omega \cup R_0$. They realized that the linear operator in this limiting problem possesses a deficiency in the resolvent estimate and belonged to a class of linear operators called \emph{almost sectorial operators} (which we will define in Section \ref{preliminares}). Nevertheless, this setback on the estimate of the resolvent did not prevent them to investigate the asymptotic behavior of the equation, which was done through a study of the global attractor for the equations.

Using this framework, Carvalho et al. in \cite{CDN} explored this deficiency in the resolvent now for the nonautonomous equation. They considered the equations \eqref{Pe*} and \eqref{P0*} in Dumbbell domains and constructed \color{black} the family of linear operators $A_{\varepsilon}(t) : D(A_{\varepsilon}(t)) \subset U_p^{\varepsilon} \rightarrow U_p^{\varepsilon}$ given by $A_{\varepsilon}(t) u = -a(t) \Delta u +u$, for $0 < \varepsilon \leq 1$, with $D(A_{\varepsilon}(t))  = \{u \in W^{2,p}(\Omega_{\varepsilon}) : \mbox{ } {\partial u }/{\partial n} = 0 \mbox{ in }\partial \Omega_{\varepsilon}\}$. For the case $\varepsilon =0$, the family of linear operators considered was $A_0(t):D(A_0(t)) \subset U_p^0 \rightarrow U_p^0$ given by $A_0(t) (w,v)  = ( -a(t)\Delta w + w, -a(t)\frac{1}{g} (gv')'+v )$ with $D(A_0(t)) = \{  (w,v) \in U_p^0 :
w  \in D( \Delta_N^{\Omega})
%\in W^{2,p}(\Omega),  \mbox{ }{\partial w}/{\partial n} = 0 \mbox{ in }\partial \Omega
;\mbox{ } (gv')'  \in L^p_g(0,1), \mbox{ }v(p_i) = w(p_i),i=0,1\}.$ Note that the domain for both operators are fixed in time.
%, $D(A_{\varepsilon}(t)) = D(A_{\varepsilon} (\cdot))$ and $D(A_0(t)) = D(A_0(\cdot))$.

Regarding the nonlinearity, they considered $F_{\varepsilon}(u) (x) = f(u(x))$, for $x \in \Omega_{\varepsilon}$ and $F_0(w,v) = (\overline{w}, \overline{v})$ where $\overline{w}(x) = f(w(x))$, if $x\in \Omega$, and $\overline{v}(x) = f(v(x))$, if $x\in R_0$.

For $\varepsilon \in (0,1]$, problem \eqref{Pe*} was written as
\begin{equation}\label{abs_e}
\begin{cases}
u_t+A_{\varepsilon}(t) u = F_{\varepsilon}(u), \quad t>\tau, \\
u(\tau) =u_0 \in U_p^{\varepsilon},
\end{cases}
\end{equation}
and, for $\varepsilon =0 $, problem \eqref{P0*} was given in its abstract form
\begin{equation}\label{abs_0}
\begin{cases}
(w,v)_t+A_{0}(t) (w,v) = F_0(w,v), \quad t>\tau, \\
(w,v)(\tau) =(w_0,v_0) \in U_p^{0}.
\end{cases}
\end{equation}

They proved that each operator of the family $\{A_0(t); t\in \R\}$ is closed, densely defined and satisfies %, for each fixed $t$, 
the estimate
\begin{equation*}
\norma{(\lambda + A_0(t))^{-1}}_{   \mathcal{L} (U_p^{0})   } \leq \dfrac{C}{|\lambda|^{\alpha}+1},
\end{equation*}
for $\lambda $ in an appropriate sector in $\C$ and $\alpha \in (0,1)$. The sector and constants being uniform in $t$. 

A family of linear operators with those properties was called \emph{uniform almost sectorial} and 
%we say that there is 
it has a deficiency in the resolvent estimate, since $\alpha <1$. This deficiency comes from the continuity condition at the junction points of the fixed part $\Omega$ with $R_0$ and it prevented the authors to use the standard theory of generation of strongly continuous semigroups and continuous linear process (provided, for instance, in \cite{pazy,sobol}).

Nevertheless, the authors in \cite{CDN} proved that each operator $-A_{0}(t)$ generates a semigroup of growth $1-\alpha$ (a concept that will be presented later in Section \ref{preliminares}) and used those semigroups to solve the linear nonautonomous evolution equation
$$\begin{cases}
(w,v)_t+A_{0}(t) (w,v) = 0, \quad t>\tau, \\
(w,v)(\tau) =(w_0,v_0) \in U_p^{{0}}.
\end{cases}$$

The solution of this problem generates a linear process of growth $1-\alpha$, denoted by $\{U_0(t,\tau):U_p^{0} \rightarrow U_p^{0}; t\geq \tau\}$, %(Definition \ref{d18.1}) 
and the function $ [\tau, \infty) \ni t \rightarrow U_0(t,\tau)(w_0,v_0)$ is not continuous at $t=\tau$ for general $(w_0,v_0)\in U_p^{0}$. In Section \ref{preliminares} we will present further details on those semigroups and process with growth.

Furthermore, the authors % in \cite{CDN} 
managed to prove that the semilinear evolution equation \eqref{abs_0} can be locally solved, even though the linear operators $A_0(t)$ presents the mentioned deficiency in the resolvent. The solution, however, possesses the same type of discontinuity at $t=\tau$ that the linear process $\{U_0 (t,\tau); t\geq \tau\}$ does.

We are interested in \color{black} continuing this study of nonautonomous reaction-diffusion equation in which features a family of almost sectorial operators. \color{black} We will compare the asymptotic dynamic of problems \eqref{abs_e} and \eqref{abs_0} and this shall be done through an analysis of attractors for both problems.

In this paper we focus on proving the existence of 
%those objects in the phase spaces $U_p^{\varepsilon}$, that is, the existence 
pullback attractor $\{ \A_{\varepsilon} (t) \subset U_p^{\varepsilon}; t\in \R  \}$ and uniform attractor $\overline{\A}_{\varepsilon} \subset U_p^{\varepsilon}$ for each $\varepsilon \in [0,1]$, with special attention to limiting case $\varepsilon = 0$, where the singular behavior appears. We will also obtain uniform estimates (in $\varepsilon$) for the attractors. 

The outline that we will follow is:
\begin{enumerate}
	\item We will first show that, under appropriate dissipativeness condition on the nonlinearity $f$, the local solution obtained in \cite{CDN} for problem \eqref{abs_0} is globally defined in time. The classical variational approach of constructing an energy functional for the equation is not possible in this case, due to the nonautonomous coefficient that multiplies the diffusion term.
	% and the fact that the  linear operator for the limiting case is not self-adjoint. 
	To overcome this problem, we will prove that the solutions of nonlinear equation are positive (that is, positive initial condition generates positive solution) and monotone (if two initial conditions satisfies an order relation, than the solutions they generate preserve this order). This is an interesting feature, considering that the solution can present discontinuity at the initial time.
	
	\item With the global existence of the solution, we will prove that the nonlinear process associated to the equation possesses a compact \color{black} set that pullback absorbs. The existence of the pullback attractor follows from the existence of such set.\color{black}
	
	\item We will prove that the compact pullback absorbing set found in (2) can be chosen uniform in $\varepsilon$ % (we will make sense of what this uniformity means, since each phase space differs from the others) 
	\color{black} and that it is also a forward attracting set. This guarantees the existence of uniform attractor for each problem and we can relate the pullback attractor and the uniform attractor obtained.
\end{enumerate}

\color{black}

Steps (2) and (3) above illustrate an interesting feature of the problem being considered, that is the fact that information of the forward dynamics can be obtained from information of the pullback dynamics. 
%As a matter of fact, we will prove the existence of a compact set that attracts bounded set for both pullback and forward dynamics. 
In general, those dynamics do not need to be related.

This paper is organized as follows: Section 2 is dedicated to the preliminaries results. We present the definitions of semigroups and process with growth $\alpha$, the definitions of almost sectorial operators and we also present the result of existence of local solution for an abstract semilinear evolution equation with almost sectorial operators. Additionally, we briefly introduce the concepts of pullback and uniform attractor, as well as results that ensure their existence. In Section 3 we begin the discussion of the equation in Dumbbell domain, and we establish the basic properties of the linear operators $A_{\varepsilon}(t)$ and nonlinearity $F_{\varepsilon}$, $\varepsilon \in [0,1]$. In Section 4 we prove the  positivity and a comparison result for the local solutions of the semilinear equation. This is done in Theorems \ref{t40.1} and \ref{t40.2}, respectively. In Section 5 we apply those results to obtain the existence of global solution (Theorem \ref{t41.1,41.2,41.4}) and the existence of pullback and uniform attractors (Theorems \ref{t41.3_41.5} and \ref{unif_attractor}).

%%%%%%%%%%%%%%%%%%%%%%%%%%%%%%%%%%%%%%

%%%%%%%%%%%%%%%%%%%%%%%%%%%%%%%%%%%%%%

%%%%%%%%%%%% PRELIMINARIES %%%%%%%%%%%%%%%

%%%%%%%%%%%%%%%%%%%%%%%%%%%%%%%%%%%%%%

%%%%%%%%%%%%%%%%%%%%%%%%%%%%%%%%%%%%%%

\section{Preliminaries}\label{preliminares}

In this section we present the concepts of \emph{semigroups of growth $\alpha$}, \emph{process of growth $\alpha$} %\emph{almost sectorial operator} 
and the result on the existence of local solution for a semilinear evolution equation with almost sectorial operator.

\subsection{Semigroups and process of growth $\alpha$: definition and properties}

\begin{definition}\label{d19.0}
	Let $X$ be a Banach space and $\alpha>0$. A family of bounded linear operators  in $X$, $\{T(t) \in \mathcal{L} (X); t > 0\}$, is a  \emph{semigroup of growth $\alpha$} if satisfies
	\begin{enumerate}
		\item $T(0) = I$.
		\item $T(t)T(s) = T(t+s)$,  for $ t,s > 0$.
		\item If $T(t)x = 0$ for all $t>0$, then $x=0$.
		\item There exists $\delta>0$ such that $\norma{t^{\alpha} T(t)}_{\mathcal{L}(X)} \leq M$, for all $0< t \leq \delta$.
		\item $X_0 = \bigcup_{t>0} T(t)X$ is dense in $X$.
	\end{enumerate}
\end{definition}

Unlike $C_0\mbox{-semigroups}$, the semigroups of growth $\alpha$ might present discontinuity at the initial time $t=0$, that is, $T(t)x \nrightarrow x$ when $t \rightarrow 0^+$ for general $x\in X$. The following properties for semigroups of growth $\alpha$ can be found in \cite{daprato,okazawa}.

\begin{lemma}\label{l2.3}
	Let $\{T(t); t>0 \}$ be a semigroup of growth $\alpha$. There exist constants $M \geq 1$ and $w\in \R$ such that, for $t>0$, $\norma{ T(t)}_{\mathcal{L}(X)} \leq M e^{wt}t^{-\alpha}$
%	\begin{equation}
%	\norma{ T(t)} \leq M e^{wt}t^{-\alpha}, \qquad \forall t>0
%	\end{equation}
	and the map $t \mapsto T(t)x$ is continuous in $(0,\infty)$.
\end{lemma}

	Let $a\in \C$ and $\varphi \in (\frac{\pi}{2}, \pi)$. We denote the sector centered in $a$ and with an angle $\varphi$ by
$$\Sigma_{\varphi, a} = \left\{ \lambda \in \C; |\arg (\lambda -a)| \leq \varphi    \right\}.$$
If $a=0$, we simply denote the sector by $\Sigma_{\varphi }$.

A linear operator  $A$ is sectorial if $-A:D(A) \subset X \rightarrow X$ is closed, densely defined linear operator, $\Sigma_{\varphi}\cup\{0\}$ belongs to its resolvent and the following estimate holds
\begin{equation*}%\label{sectorial_limitation}
\norma{(\lambda +A)^{-1}}_{\mathcal{L}(X)} \leq \dfrac{M}{|\lambda|}, \quad \lambda \in \Sigma_{\varphi}.
\end{equation*}

The operator $-A$ generates a $C_0\mbox{-semigroup}$ which is given by the following integral
\begin{equation*}%\label{semigrupo_analitico}
T_{-A}(t) = \dfrac{1}{2\pi i} \ds\int_{\Gamma} e^{\lambda t}(\lambda +A)^{-1} d\lambda,
\end{equation*}
where $\Gamma$ is the boundary of the sector $\Sigma_{\varphi} $. The above characterization of the semigroup allows one to prove that the semigroup is analytic and it regularizes as the time evolves (see \cite{pazy} for more details).

However, in some applications, as the equation in Dumbbell domains that we consider here, the estimate on the resolvent of the operator $-A$ possesses a deficiency. In this case, $A$ will be an \emph{almost sectorial operator}, which we define next.

\begin{definition}\label{d12.1}
	A linear operator $A:D(A) \subset X \rightarrow X$ is \emph{almost sectorial} if 
	\begin{enumerate}
		\item $-A$ is closed and densely defined.
		\item There exist $ \varphi \in \left(  \frac{\pi}{2}, \pi  \right) $ and $M>0$ such that $\Sigma_{\varphi} \cup \{ 0 \} \subset \rho(-A)$
		and, for all $ \lambda \in \Sigma_{\varphi}$,
		\begin{equation*}
		\norma{(\lambda +A)^{-1}}_{\mathcal{L}(X)} \leq \dfrac{M}{|\lambda|^{\alpha}},
		\end{equation*}
		for some $\alpha \in (0,1)$.
	\end{enumerate}
\end{definition}

Despite the fact that $\alpha$ is strictly less than $1$, we can still define for an almost sectorial operators the family
\begin{equation*}%\label{almost_sectorial}
T_{-A}(t) = \dfrac{1}{2\pi i} \ds\int_{\Gamma} e^{\lambda t}(\lambda +A)^{-1} d\lambda.
\end{equation*}

In \cite{DDDII} the authors proved that the above integral converges and this family of linear operators defines a semigroup of growth  $1-\alpha$. Such semigroup is continuous in $(0,\infty)$ and, following the same ideas used in \cite{pazy}, this semigroup is analytic in any interval $[\delta, \infty)$, for $\delta >0$. For further details and examples of semigroups of growth $\alpha$, we recommend \cite{daprato, krein}. In \cite{okazawa} a result on the generation of those type of semigroups is given and in \cite{PS} a functional calculus for almost sectorial operators is developed.

When dealing with nonautonomous evolution equations with almost sectorial operators, as we shall do in the sequence, we will also work with a 2-parameters family of linear operators with growth $ \alpha$, which we define next.

\begin{definition}\label{d18.1}
	Let $X$ be a Banach space and $\alpha >0$. A family $\{U(t,s) \in \mathcal{L}(X);t> s\}$  is a \emph{process of growth $\alpha$} if 
	
	\begin{enumerate}
		\item $U(t,t) = Id$.
		\item $U(t, \tau)U(\tau, s) = U(t,\tau) \circ U(\tau ,s) = U(t,s)$, for all $s<\tau < t$.
		\item There exists $M>0$ such that $\norma{(t-s)^{\alpha}U(t,s)}_{\mathcal{L}(X)} \leq M$, for all $ t > s$.
		\item $(t,s,x) \rightarrow U(t,s)x$ is continuous for $t>s$ and for all $ x \in X$. 
	\end{enumerate}
\end{definition}

	As it happens for the semigroups of growth $\alpha$, the process of growth $\alpha$ might be discontinuous at the initial time $t= \tau$.

%%%%%%%%%%%%%%%%%%%%%%%%%%%%%%%%%%%%%%

%%%%%%%%%%%%%%%%%%%%%%%%%%%%%%%%%%%%%%

\subsection{Semilinear evolution equations with almost sectorial operators: Abstract theory}

Consider the semilinear evolution equation
\begin{equation}\label{SEE}
\begin{cases}
u_t + A(t)u = F(u), \quad t>\tau, \\ 
u(\tau) = u_0\in X,
\end{cases}
\end{equation}
where $\{A(t): D(A(t)) \subset X \rightarrow X, \mbox{ } t\in \R\}$ is a family of  \emph{uniformly almost sectorial operators}, that is, for all $t\in \R$
\begin{enumerate}
	\item $-A(t): D(A(t)) \subset X \rightarrow X$ is closed, densely defined and $D(A(t)) = D = X^1$, $\forall t \in \R$.
	
	\item There exists $\varphi \in \left( \frac{\pi}{2}, \pi \right)$ and constants $C>0$, $\alpha \in (0,1)$ such that $\Sigma_{\varphi }\cup \{0\} \subset \rho (-A(t)),$
	and 
	\begin{equation}\label{13.I}
	\norma{(\lambda + A(t))^{-1}}_{\mathcal{L}(X)} 
	\leq
	 \dfrac{C}{|\lambda|^{\alpha}},\quad \forall \lambda \in \Sigma_{\varphi}.
	\end{equation}
	\item The constants $\varphi$, $\alpha$ and $C$ are the same for the entire family $\{A(t); t\in \R\}$.
\end{enumerate}

We also assume that $\{A(t); t\in \R\}$ is \emph{uniformly H\"{o}lder continuous}, that is, there exist $\delta \in (0,1]$ and a constant $C>0$ such that 
\begin{equation*}\label{HC}
\norma{[A(t)-A(\tau)] A(s)^{-1}}_{\mathcal{L}(X)} \leq C|t-\tau|^{\delta}, \quad \forall t,\tau,s \in \R.
\end{equation*} 

By taking $\tau = s$ we conclude that
%$$\norma{ A(t)A(s)^{-1} - I} \leq C (t-s)^{\varepsilon}$$
\begin{equation*}
\norma{A(t)A(s)^{-1} }_{\mathcal{L}(X)}  \leq 1+C(t-s)^{\delta}
\end{equation*}
and that $A(t)A(s)^{-1}$ is a bounded linear operator in $X$. From the fact that $0 \in \rho (-A(t))$ and from the continuity of the function $\rho{(-A(t))} \ni \lambda \mapsto (\lambda +A(t))^{-1}$, the inequality \eqref{13.I} is equivalent to 
\begin{equation*}%\label{13.II}
\norma{(\lambda + A(t))^{-1} }_{\mathcal{L}(X)} \leq \dfrac{C}{1+|\lambda|^{\alpha}}, \quad \forall \lambda\in\Sigma_{\varphi}\cup\{0\}.
\end{equation*}

Still from \eqref{13.I} and from the resolvent equality, we have 
\begin{equation*}%\label{13.III}
\norma{A(t)(\lambda + A(t))^{-1}}_{\mathcal{L}(X)} \leq 1+C|\lambda|^{1-\alpha}, \quad \forall \lambda\in\Sigma_{\varphi}\cup\{0\}.
\end{equation*}

Under those conditions, for a fixed $\tau \in \R$, the operator $-A(\tau)$ generates a semigroup of growth $1-\alpha$ given by
\begin{equation}\label{semigrupo}
T_{-A(\tau)}(t)x = \dfrac{1}{2\pi i} \int_{\Gamma} e^{\lambda t} (\lambda + A(\tau))^{-1} x d\lambda,
\end{equation}
where $\Gamma$ is the boundary of the sector $\Sigma_{\varphi}$ and $\norma{T_{-A(\tau)}(t)} \leq C t^{\alpha-1}$, for $t>0$.
%\begin{equation}\label{13.IV}
%\norma{t^{1-\alpha}T_{-A(\tau)}(t)} \leq C, \quad \forall t>0.
%\end{equation}

Since $ 0 \in \rho (-A(\cdot))$ and the resolvent is an open set, there exists $\xi >0$ such that the family $\{ -\xi I+ A(t); t \in \R\}$ is uniformly almost sectorial. Therefore, as it is done in \cite{pazy} for sectorial operators, we conclude that the semigroup has an exponential decay
\begin{equation*}%\label{exp_decay}	
\norma{T_{-A(\tau)}(t) }_{\mathcal{L}(X)} \leq C t^{\alpha-1} e^{-\xi t}, \quad \forall t>0.
\end{equation*}

	In \cite{CDN}, to obtain the existence of local solution for \eqref{SEE}, the authors studied the problems in three stages: first the autonomous linear case, then the nonautonomous linear equation and, at last, the semilinear case. Next, we briefly mention this path followed by the authors.

%%%%%%%%%%%%%%%%%%%%%%%%%%%%%%%%%%%%%%

%%%%%%%%%%%%%%%%%%%%%%%%%%%%%%%%%%%%%%

\subsection{The autonomous linear problem}\label{PLH}

Consider the autonomous initial value problem
\begin{equation}\label{13.XIII}
\begin{cases} u_t+ A(\tau) u = 0, \quad t> 0, \\ u(0) = u_0 \in X,
\end{cases}
\end{equation}
where $\tau$ is fixed. The next result can be found in \cite{DDDII} and it claims that the semigroup generated by $-A(\tau)$ gives a solution for the problem.

\begin{lemma}\label{l13.7}
	Let $\{T_{-A(\tau)}(t); t  > 0\}$ be the semigroup of growth $1-\alpha$ defined above. Then $T_{-A(\tau)}(t):(0,\infty) \rightarrow \mathcal{L}(X)$ is differentiable and $\frac{d}{dt}T_{-A(\tau)}(t) = -A(\tau) T_{-A(\tau)}(t),$
	that is, for all $u_0 \in X$,
	$$\dfrac{d}{dt}T_{-A(\tau)}(t)u_0 + A(\tau) T_{-A(\tau)}(t) u_0 = 0, \quad \forall t>0,$$
	and $u(t) = T_{-A(\tau)}(t) u_0$ is a solution of \eqref{13.XIII}.
\end{lemma}

We already mentioned that semigroups of growth $1-\alpha$  are not necessarily continuous at $t=0$. However, it was proved in \cite{DDDII} that, for initial conditions in $D = D(A(\cdot))$, the continuity at $t=0$ follows.

\begin{lemma}\label{l2.6}
	If $u_0 \in D = D(A(\cdot))$ then 
	\begin{enumerate}
		\item $\norma{T_{-A(\tau)}(t)u_0 - u_0 }_{X} \rightarrow 0$ when $t \rightarrow 0^+$.
		\item $T_{-A(\tau)}(t) A(\tau) u_0 = A(\tau) T_{-A(\tau)}(t) u_0$, $\forall t>0$. In this case, the solution $u(t) = 	 T_{-A(\tau)}(t) u_0$ is continuously differentiable in $(0,\infty)$.
	\end{enumerate}
	
\end{lemma}

%%%%%%%%%%%%%%%%%%%%%%%%%%%%%%%%%%%%%%

%%%%%%%%%%%%%%%%%%%%%%%%%%%%%%%%%%%%%%

\subsection{The nonautonomous linear problem}\label{PLNH}
Next, we consider 
\begin{equation}\label{14.II}
\begin{cases}
u_t+A(t)u = 0, \quad t> \tau,
 \\ u(\tau) = u_0 \in X, \quad \tau \in \R.
\end{cases}
\end{equation}

In \cite{CDN}, assuming $\alpha + \delta >1$, where $\delta >0$ comes from the H\"{o}lder continuity \eqref{HC}, the authors proved the existence of a family of solution operators $\{U(t, \tau) \in \mathcal{L}(X); \mbox{ } t\geq \tau\}$ , given by
\begin{equation}\label{14.IV}
U(t,\tau) = T_{-A(\tau)}(t-\tau) + \int_{\tau}^t U(t,s) [ A(\tau) - A(s)]T_{-A(\tau)}(s-\tau)ds.
\end{equation}

Such family is a linear process of growth $1-\alpha$, that is,
\begin{equation*}
\norma{U(t,\tau)}_{\mathcal{L}(X)} \leq C(t-\tau)^{\alpha-1} , \quad \forall t> \tau.
\end{equation*}

The authors of \cite{CDN} also exhibited another integral formulation for this process that makes it easier to deal with %the analisys of the differentiability of solutions given by the process. 
the analysis of some properties of the process. Let $\varphi_1$ and $\Phi$ be
\begin{equation*}%\label{31.I}
\varphi_1(t,\tau) = [A(\tau) - A(t)] T_{-A(\tau)}(t-\tau)
\end{equation*}
\begin{equation*}%\label{31.II}
\Phi (t,\tau)= \varphi_1 (t,\tau) + \ds\int_{\tau}^t \Phi (t,s) \varphi_1(s,\tau) ds = \varphi_1(t,\tau) + \ds\int_{\tau}^t \varphi_1(t,s) \Phi (s,\tau)ds.
\end{equation*}

Then, if $\alpha +\frac{ \delta}{2}>1$, the linear process \eqref{14.IV} can also be given by
\begin{equation}\label{31.III}
U(t,\tau) = T_{-A(\tau)}(t-\tau)+ \ds\int_{\tau}^t T_{-A(s)}(t-s) \Phi (s,\tau) ds.
\end{equation}

%Vale ressaltar que em \cite{CDN} apenas a primeira igualdade em \ref{31.II} está demonstrada. A segunda igualdade segue de forma idêntica à feita na página 28 de \cite{Tese_Marcelo}. 

The operators $\varphi_1(t,\tau)$ and $\Phi(t,\tau)$ are continuous in $\{ (t,\tau) \in \R^2; \mbox{ }t> \tau\}$ and satisfy 
\begin{equation}\label{31.IV}\norma{\varphi_1(t,\tau)}_{\mathcal{L}(X)} \leq C(t-\tau)^{\alpha+\delta -2},%e^{-\xi t}
\end{equation} \begin{equation}\label{31.V}\norma{\Phi(t,\tau)}_{\mathcal{L}(X)} \leq C(t-\tau)^{\alpha+\delta -2}.
\end{equation}

As it happens for the semigroup, if $x \in D(A(\cdot))$ then the process is continuous at the initial time.

\begin{proposition}\label{p31.12}
	If $x\in D(A(\cdot))$, then $U(t,\tau)x \xrightarrow{t\rightarrow \tau^{+}} x.$
	
	\begin{proof} Note that 
		$$\norma{U(t,\tau)x-x}_{X} \leq \norma{T_{-A(\tau)}(t-\tau)x-x}_{X} + \norma{ \ds\int_{\tau}^t T_{-A(s)}(t-s) \Phi (s,\tau)xds }_{X}.$$
		
		The first term at the right side goes to zero as a consequence of Lemma \ref{l2.6} and, using \eqref{31.IV} and \eqref{31.V} for the second term, we have
		\begin{align*}
		\norma{ \ds\int_{\tau}^t T_{-A(s)}(t-s) \Phi (s,\tau)xds }_{\mathcal{L}(X)} & \leq C \ds\int_{\tau}^t (t-s)^{\alpha-1} (s-\tau)^{\alpha+\delta-2} ds \\
		& \leq C \mathcal{B}( \alpha , \alpha+\delta-1) (t-\tau)^{2\alpha + \delta -2} \xrightarrow{t\rightarrow \tau^{+}} 0, 
		\end{align*}
		where $\mathcal{B} (\cdot, \cdot): (0,\infty) \times (0,\infty)
		\rightarrow \R$ is the \emph{Beta function} 
		$$ \mathcal{B} (a,b) = \int_0^{1} t^{a-1}(1-t)^{b-1}dt$$
		and $2\alpha + \delta -2 = 2 \left( \alpha +\frac{\delta}{2}-1 \right) >0$.	\end{proof}
\end{proposition}

%%%%%%%%%%%%%%%%%%%%%%%%%%%%%%%%%%%%%%

%%%%%%%%%%%%%%%%%%%%%%%%%%%%%%%%%%%%%%

\subsection{Existence of local solution for the semilinear problem}\label{sol_loc} 

Finally we consider the semilinear case 
\begin{equation}\label{18.II}
\begin{cases}
u_t+A(t)u = F(u), \quad t> \tau,\\
u(\tau) = u_0 \in X,
\end{cases}
\end{equation}
where the family $\{A(t);\mbox{ } t\in \R\}$ is uniformly almost sectorial and uniformly H\"{o}lder continuous.

We also assume that there exists another Banach space $Y$ and $0 < \beta <1$ such %that the family $\{A(t);t\in \R\}$ has its image in $Y$, 
that $\beta + \delta >1$ and the resolvent satisfies the following estimate
\begin{equation}\label{25.*}
\norma{ (\lambda+A(t))^{-1}}_{\mathcal{L}(Y,X)} \leq \dfrac{C}{|\lambda|^{\beta}+1},\qquad \forall \lambda \in \Sigma_{\varphi}\cup\{0\}.
\end{equation}

Using \eqref{semigrupo}, \eqref{14.IV}, \eqref{25.*} and the fact that $0 \in \rho (-A(\tau))$, we have, for $t>0$,
\begin{equation*}%\label{25.I}
\begin{matrix}
\norma{T_{-A(\tau)}(t)}_{\mathcal{L}(Y,X)} \leq Ct^{\beta-1}	& \mbox{ and } & \norma{T_{-A(\tau)}(t)}_{\mathcal{L}(Y,X)} \leq Ct^{\beta-1}e^{-\xi t}
\end{matrix}
\end{equation*}
and, for $t>\tau$,
	\begin{equation*}
	\norma{U(t,\tau)}_{\mathcal{L}(Y,X)} \leq C (t-\tau)^{\beta-1}.
	\end{equation*}

We assume that the nonlinearity $F$ satisfies: $F:X\rightarrow Y$ and there exist constants $ C>0 \mbox{ and }\rho \geq 1$
 such that, for all $u,v \in X$,
\begin{equation}\label{25.V}
\norma{F(u) - F(v) }_{Y} \leq C \norma{u-v}_{X}\left( 1+\norma{u}_{X}^{\rho-1}+\norma{v}_{X}^{\rho-1}\right), 
\end{equation}
\begin{equation}\label{25.VI}
\norma{F(u)}_{Y} \leq C (1+ \norma{u}_{X}^{\rho}). 
\end{equation}
The exponents $\alpha, \beta, \delta$ and $\rho$ must satisfy
\begin{equation}\label{25.VII}
\dfrac{1}{2}< \alpha, \beta <1,
\end{equation}
%\begin{equation}\label{25.VIII}
%\alpha, \beta > \dfrac{1}{2}
%\end{equation}
\begin{equation}\label{25.IX}
\dfrac{\beta}{1-\alpha}>1 \quad \mbox{ and }\quad 1 \leq \rho < \dfrac{\beta}{1-\alpha},
\end{equation}
\begin{equation}\label{25.X}
\alpha + \delta > 1 \quad \mbox{ and }\quad\beta + \delta >1.
\end{equation}

\begin{definition}\label{d25.3}
	A function $u:(\tau , \tau+t_0) \rightarrow X$ is a \emph{mild solution} for \eqref{18.II} if 
	$$(\tau, \tau+t_0) \ni t \mapsto w(t):= u(t) - U(t,\tau)u_0 \in X$$ is continuous, $\lim_{t \rightarrow \tau^+} \norma{w(t)}_{X} = 0$ and
	$$w(t) = u(t)- U(t,\tau)u_0 = \int_{\tau}^t U(t,s)F(u(s)) ds, \quad \forall t \in (\tau, \tau+t_0).$$
	
\end{definition}

We point out here the main idea used in proving the existence of local solution, which is the definition of a contraction map in an appropriate space. This idea will be important in the sequence. For a complete proof, the reader may consult \cite[pages 34-36]{CDN} or \cite[pages 184-187]{DDDII} (this second one deals with the autonomous case, but the general setting and the ideas applied there are essentially the same).

Since $\{U(t,\tau);  t> \tau \}$, solution of \eqref{14.II}, is discontinuous at $t = \tau$ we expect that the same happens for solutions of the semilinear case \eqref{18.II}. Hence, given $\tau \in\R$ and $u_0 \in X$, we will search for mild solutions in the following space 
% if $\tau \in \R$ and $u_0 \in X$, we search for $t_0>0$ such that exists $u: (\tau,\tau+t_0) \mapsto X$ mild solution \ref{18.II}. Para isso, consideremos o seguinte espa\c{c}o
$$K(t_0,u_0) = \left\{ v \in \mathcal{C}((\tau, \tau+t_0), X); \mbox{ } \sup_{t \in (\tau,\tau+t_0)} \norma{v(t)-U(t,\tau)u_0}_{X} \leq \mu\right\},$$
where $\mu>0$ and $t_0>0$ will be chosen later.
Such space is a Banach space with the norm
$\norma{\phi }_{K(t_0,u_0)} = \sup_{t \in (\tau,\tau+t_0)} \norma{\phi(t)-U(t,\tau)u_0}_{X}.$

\begin{theorem}\cite[Theorem 3.1]{CDN}\label{t25.4}
	Suppose \eqref{25.V} up to \eqref{25.X} are satisfied. Then, for every $u_0 \in X$, there exists $t_0>0$ small enough so that the initial value problem
	\begin{equation*}
	\begin{cases}
	u_t+A(t)u = F(u), \quad t> \tau,\\
	u(\tau) = u_0 \in X,
	\end{cases}
	\end{equation*}
	has a mild solution in $K(t_0,u_0)$. Furthermore, we can extend this mild solution to a maximal interval $(\tau, \tau+t_{\max}(u_0))$. Both $t_0$ and $t_{\max}(u_0)$ depend on $u_0$, but can be chosen uniformly for $u_0 $ in bounded sets.
\end{theorem}

To prove this theorem, the authors considered the operator
$T: K(t_0,u_0) \rightarrow \mathcal{C}((\tau, \tau+t_0),X)$ defined by
$$(Tv)(t) : = U(t,\tau)u_0 + \ds\int_{\tau}^t U(t,s)F(u(s))ds.$$
For $t_0>0$ small enough, they proved that this map is a contraction and its fixed point is a mild solution of \eqref{18.II}. Also, for any $v \in K(t_0,u_0)$,
{\small \begin{equation*}%\label{25.XI}
	\begin{split}
	(t-\tau)^{1-\alpha} \norma{v(t)}_{X} & \leq (t-\tau)^{ 1-\alpha } \norma{v(t)- U(t,\tau)u_0}_X + (t-\tau)^{1-\alpha} \norma{U(t,\tau)v_0}_X \\
	& \leq (t-\tau)^{1-\alpha} \mu + C \norma{u_0}_X  \leq t_0^{1-\alpha} \mu + C \sup_{u_0 \in B} \norma{u_0} \leq k,
	\end{split}
	\end{equation*}}where $B \subset X$ is a bounded set that contains $u_0$. It follows that $t_0$ can be taken uniformly for initial conditions in bounded set, $u_0 \in B$.

Using the same ideas that appears in \cite{DDDII} for the autonomous case, we managed to prove the following proposition.

\begin{proposition}\label{t25.5}
	In the conditions of Theorem \ref{t25.4}, for each $u_0 \in X$ there exists a unique mild solution of \eqref{18.II}, $u:(\tau, \tau+t_{max}(u_0)) \rightarrow X$, defined in a maximal interval, which satisfies
	
	\begin{enumerate}
		\item $t_{max}(u_0) = +\infty$ or $ \limsup_{t \rightarrow\tau+t_{max}(u_0)} \norma{u(t)}_{X} = +\infty.$
		
		\item The solution depends continuously on the initial data in the following sense: if $u_0 \in X$ and $t^*< t_{max}(u_0)$, then there exists $\delta >0$ small enough such that 
		$$\norma{u_0-v_0}_{X} < \delta \mbox{ and } t\in (\tau, \tau+t^*) \Rightarrow \norma{u(t,\tau,u_0) - v(t,\tau,v_0) }_{X} \leq C(t-\tau)^{\alpha-1} \norma{u_0-v_0}_{X},$$
		that is, the solutions of \eqref{18.II} behave as the solutions of the associated linear problem.
	\end{enumerate}

\end{proposition}

%%%%%%%%%%%%%%%%%%%%%%%%%%%%%%%%%%%%%%

%%%%%%%%%%%%%%%%%%%%%%%%%%%%%%%%%%%%%%

\subsection{Pullback and uniform attractors}\label{attractors} 

Let $X$ be a Banach space and $\{S(t,\tau) :X \rightarrow X; t\geq s\}$ a family of operators satisfying:
\begin{enumerate}
	\item $S(t,t) = I_X$, for all $t \in \R$.
	\item $S(t,s) = S(t, \tau) S(\tau, s)$, for all $t\geq \tau \geq s$, $s\in \R$.
	\item $(s,\infty) \ni t \mapsto S(t,s)x$ is continuous for all $x \in X$.
\end{enumerate}
Such family is called an \emph{evolution process in $X$} and we also denote it by $S(\cdot,\cdot)$.

 Given any $x\in X$, 
 %let $S(t,\tau)x$ be the evolution of the system that started at the inicial condition $x$ at time $s$ and evolved $t-s$ in time. %and the final state $S(t,s)x$
there are two distinct manners of studying the asymptotic dynamics of such evolution process: One called the \emph{pullback dynamics} that basically fixes the final time $t$ and evaluate what happens to $S(t,s)x$ when
%$s$ is each time more negative, that is,
 $s\rightarrow -\infty$ and the other called \emph{forward dynamics}, which consider $S(t,s)x$ when $s$ is fixed and $t \rightarrow \infty$.

 %
 %Pullback
 %

 The pullback dynamics can be described by an object in the phase space called \emph{pullback attractor}. We recall in the sequence some basic concepts and results of the theory of pullback attractor. We refer to \cite{Carvalho} and references therein for further details. 
 
 Furthermore, throughout the text, we shall use the Hausdorff semidistance to compare the distance between two sets in the phase space $X$, that is,  given $A,B \subset X$, the \emph{Hausdorff semidistance} between $A$ and $B$ is given by
 $$dist(A,B) = \sup_{a \in A} \inf_{b \in B} d(a,b).$$
 
 \begin{definition}
 	Let $S(\cdot, \cdot)$ be a process. A family $A(\cdot)= \{A(t)\subset X; t\in \R\}$ \emph{pullback attracts $B \subset X$} %under the process $S(\cdot, \cdot)$
 	 if, for each $t \in \R$,  
 		$dist(S(t,s)B, A(t)) \rightarrow 0 \mbox{ when } s\rightarrow - \infty.$
 \end{definition}
 
 \begin{definition}
 	The \emph{pullback attractor} of $S(\cdot, \cdot)$ is a family $\A(\cdot) = \{\A(t) \subset X; \mbox{ } t \in \R\}$ that satisfies:
 	
 	\begin{enumerate}
 		\item $\A(t)$ is compact for all $t\in \R$.
 		\item $\A(\cdot)$ is \emph{invariant} by $S(\cdot, \cdot)$, that is, $S(t,s)\A(s) = \A(t)$, for all $t\geq s$, $s\in \R$.
 		\item $\A(\cdot)$ pullback attracts bounded sets of $X$.
 		\item $\A(\cdot)$ is the minimal closed family that satisfies $(3)$.
 	\end{enumerate}
 \end{definition}

\begin{theorem}(See \cite{Carvalho})\label{t21.1}
	Let $S(\cdot,\cdot)$ be a process. The statements below are equivalent:
	
	\begin{enumerate}
		\item $S(\cdot,\cdot)$ has a pullback attractor $\A(\cdot)$.
		
		\item There exists a family of compact sets $K(\cdot)$ that pullback attracts bounded sets of $X$.
	\end{enumerate}

\end{theorem}

 When it comes to the forward dynamics, the object that describes the asymptotic behavior is called \emph{uniform attractor}. We present some concepts and results next, referring to \cite{B-V,Carvalho,ChepyzhovVishik,WangLiQin} for more details.

 \begin{definition}
 	Let $S(\cdot, \cdot)$ be a process. A set $A \subset X$ is  \emph{uniformly attracting} if, for any  $B \subset X$ bounded,
 	$$ \lim_{t \rightarrow \infty}  \left(  \sup_{s \in \R} dist \left(  S(t + s , s )  B ,  A \right)\right) =0.$$
 \end{definition}

 \begin{definition}
 The minimal compact attracting set for $S(\cdot, \cdot)$ is  called \emph{uniform attractor}, denoted by $\overline{\A}$.% Note that the uniform attractor is a fixed set in the phase space and the attraction happens forward in this case.
 \end{definition}

\begin{theorem}(See \cite{Carvalho})\label{exist_uni_atr}
	Let $S(\cdot,\cdot)$ be a process. The statements below are equivalent:
	
	\begin{enumerate}
		\item $S(\cdot,\cdot)$ has uniform attractor $\overline{\A}$.
		
		\item There exists a compact uniformly attracting set $K$ for the process $S(\cdot,\cdot)$.
	\end{enumerate}
	
	In those cases, $S(\cdot,\cdot)$ also has a pullback attractor $\A(\cdot)$ and 
	\[\cup_{t\in \R}   \A(t) \subset \overline{\A}.  \]
	
\end{theorem}

 \color{black}
 
 \begin{remark}
 	The results on existence of attractors provided in \cite{Carvalho} are all obtained under the assumption that  $[s,\infty) \ni t \mapsto S(t,s)x$ is continuous for all $x \in X$, which includes left extreme of the interval, $t=s$. The processes that we will obtain in this work will not be continuous at the initial time. Nevertheless, the results on characterization and existence of attractors remain valid and no significant changes in the proofs of those results are necessary. %for the proofs of such results.
 \end{remark}

 %%%%%%%%%%%%%%%%%%%%%%%%%%%%%%%%%%%%%%
 
 %%%%%%%%%%%%%%%%%%%%%%%%%%%%%%%%%%%%%%
 
 %%%%%%%%%%%%% SECTION 3 %%%%%%%%%%%%%%%%%
 
 %%%%%%%%%%%%%%%%%%%%%%%%%%%%%%%%%%%%%%
 
 %%%%%%%%%%%%%%%%%%%%%%%%%%%%%%%%%%%%%%

\section{Nonautonomous reaction-diffusion equation in Dumbbell Domains}\label{Aplication}

Consider the equations \eqref{Pe*} and \eqref{P0*} presented in the Introduction, each one taking place in a different domain $\Omega_{\varepsilon}$. We require the following properties on those sets that will determine the shape of the channels $R_{\varepsilon}$ and the way they collapse at the line segment $R_0$:

\begin{enumerate}
	\item There exists $l>0$ such that 
	$$\Omega \cap \{0<s<1; \mbox{ } |x'|< l\} = \emptyset$$
	$$\{(s,x'); \mbox{ } s^2+|x'|^2 < l^2 \mbox{ and } s<0\} \subset \Omega$$
	$$\{(s,x'); \mbox{ } (s-1)^2+|x'|^2 < l^2 \mbox{ and } s>1\} \subset \Omega$$
	$$ \{(0,x'); \mbox{ }|x'|<l\} \cup \{(1,x'); \mbox{ }|x'|<l\} \subset \partial \Omega.$$
	
	\item For $\varepsilon =1$ we have $R_1 =\{ (s,x'); \mbox{ } 0\geq s \geq 1 \mbox{ and }x'\in \Gamma_s^1\},$
	where $\Gamma_s^1 \subset R^{N-1}$ is $\mathcal{C}^1\mbox{-diffeomorphic}$ to the unitary ball of $\R^{N-1}$, that is, for every $ s \in [0,1]$, there exists a diffeomorphism of class $\mathcal{C}^1$, $L_s:B(0,1) \rightarrow \Gamma_s^1.$
	Besides, we assume that
	\begin{align*}
	L:(0,1) \times B(0,1) &\rightarrow \quad R_1 \\
	(s,x')\quad \quad  & \mapsto (s,L_s(x'))
	\end{align*}
	is a $\mathcal{C}^1\mbox{-diffeomorphism}$.
	
	\item Still for $\varepsilon =1 $, we will denote by $g(s) = |\Gamma_s^1| $ the $(N-1)\mbox{-Lebesgue measure}$ of the subset $\Gamma_s^1$. Due to the smoothness of $R_1$, we assume that $g$
	is smooth in the interval $(0,1)$ and that there exists $d_0,d_1 >0$ such that $d_0 \leq g(s) \leq d_1$, for all $s\in [0,1].$ This $g$ is the function that features in the limit equation \eqref{P0*} and the weight that appears in the norm of $U_p^{0}$.
	
	\item For $\varepsilon \in (0,1)$, we define $R_{\varepsilon} = \{(s,\varepsilon x'); \mbox{ } (s,x') \in R_1\}$ and this implies that $R_{\varepsilon}$ collapses into $R_0 = \{(s,0); \mbox{ } s\in [0,1]\}.$
\end{enumerate}

The phase spaces in which we shall consider the abstract equations \eqref{abs_e} and \eqref{abs_0} are $U_p^{\varepsilon}$ and $U_p^{0}$, defined in the Introduction, with norms \eqref{norma_e} and \eqref{norma_0}, respectively.

Note that if $u \in U_p^{\varepsilon}$ is such that $u(s,y)$ does not depend on $y$ when $(s,y) \in R_{\varepsilon}$, then
{\small 	\begin{align*}
	\int_{R_{\varepsilon}} |u(s,y)|^p dsdy & = \int_0^1 \int_{\Gamma_s^{\varepsilon}} |u(s,y)|^p dy ds  = \int_0^1 \int_{\Gamma_s^1} |u(s,\varepsilon x')|^p \varepsilon^{N-1} dx'ds = \varepsilon^{N-1} \int_0^1 g(s) |u(s)|^p ds,
	\end{align*}}
that is, $\norma{u}_{U_p^{\varepsilon}}^p = \norma{u}_{L^p(\Omega)}^p + \int_0^1 g(s) |u(s)|^p ds.$

%%%%%%%%%%%%%%%%%%%%%%%%%%%%%%%%%%%%%%

%%%%%%%%%%%%%%%%%%%%%%%%%%%%%%%%%%%%%%

\subsection{Existence of local mild solution for the problem}

To show the existence of local mild solution for problems \eqref{abs_e} and \eqref{abs_0}, we suppose that the functions $a: \R \rightarrow [c_0,c_1]$, with $0<c_0< c_1$, and $f: \R \rightarrow \R$ %
 satisfy the conditions \eqref{a_holder} and \eqref{G**}. The estimate on the growth of the function $f$ implies that, for $t,s \in \R$, 
\begin{equation}\label{26.*}
|f(s)-f(t)| \leq C|t-s|(1+|t|^{\rho -1}+|s|^{\rho -1}) \mbox{ and }|f(t)| \leq C (1+|t|^{\rho}),\\
\end{equation}
where $\rho \geq 1$.

Under these conditions, the family of operators $\{A_0(t); t \in \R\}$ defined in the Introduction is uniformly almost sectorial. The following proposition is stated in \cite{CDN} and the idea of its proof (for the autonomous case) can be seen in \cite{DDDII}.

\begin{proposition}\label{p22.1}
	The family $\{A_0(t); t\in \R\}$ satisfies:
	\begin{enumerate}
		\item $A_0(t)$ is closed and densely defined.
		\item $A_0(t)$ has a compact resolvent.
		\item There exist $\varphi \in \left( \frac{\pi}{2}, \pi \right)$ and $C>0$ (independent of $t$) such that $\Sigma_{\varphi} \subset \rho (-A_0(t))$ and, for $\frac{N}{2}<q \leq p$,
		\begin{equation*}%\label{22.I}
		\norma{(\lambda +A_0(t))^{-1}}_{\mathcal{L}(U_q^{0},U_p^{0})} \leq \dfrac{C}{|\lambda|^{\alpha}+1},
		\end{equation*}
		\begin{equation*}%\label{22.II}
		\norma{A_0(t) (\lambda +A_0(t))^{-1}}_{\mathcal{L}(U_p^0)} \leq C(1+|\lambda |^{1-\tilde{\alpha}}	)	,
		\end{equation*}
		for each $0<\alpha <1-\frac{N}{2q}-\frac{1}{2}\left( \frac{1}{q} - \frac{1}{p}\right)$, $0< \tilde{\alpha}<1-\frac{N}{2p}<1$ and $\lambda \in \Sigma_{\varphi}$.
		\item If $B_0(t)$ is the realization of $A_0(t)$ em $\mathcal{C}( \overline{\Omega}) \oplus \mathcal{C}(0,1)$, then $B_0(t)$ is a sectorial operator in $\mathcal{C}( \overline{\Omega}) \oplus \mathcal{C}(0,1)$ with compact resolvent. Therefore, $-B_0(t)$ generates a $C_0-$semigroup analytic in $\mathcal{C}( \overline{\Omega}) \oplus \mathcal{C}(0,1)$.
	\end{enumerate}
	
\end{proposition}

Besides that, it follows from the H\"{o}lder continuity of the function $a:\R \rightarrow [c_0,c_1]$ that the family $\{ A_0(t); t\in \R \}$ is also H\"{o}lder continuous 
$$\norma{[A_0(t)-A_0(\tau)] A_0(s)^{-1}}_{\mathcal{L}(U_p^0)} \leq M (t-\tau)^{\delta}, \qquad \forall t, \tau, s \in \R.$$

If $\rho = \frac{p}{q}$, inequalities in \eqref{26.*} imply %the operator $F_0:U_p^0 \rightarrow U_q^0$ satisfies
\begin{align*}
\norma{F_0(w,v)  - F_0(\tilde{w},\tilde{v})}_{U_q^0} 
%& \leq C (\norma{ w -\tilde{w}}_{L^p(\Omega)} + \norma{v-\tilde{v}}_{L^p_g(0,1)})(1+ \norma{w}_{L^p(\Omega)}^{\rho-1} + \norma{v}_{L^p_g(0,1)}^{\rho-1} +\norma{\tilde{w}}_{L^p(\Omega)}^{\rho-1} + \norma{\tilde{v}}_{L^p_g(0,1)}^{\rho-1}) \\
%& \leq C (\norma{ w -\tilde{w}}_{L^p(\Omega)} + \norma{v-\tilde{v}}_{L^p_g(0,1)})(1+ \norma{(w,v)}_{U_p^0}^{\rho-1} +\norma{ (\tilde{w},\tilde{v}) }_{U_p^0}^{\rho-1})\\
\leq C \norma{ (w,v) - (\tilde{w},\tilde{v})}_{U_p^0}(1+ \norma{(w,v)}_{U_p^0}^{\rho-1} +\norma{ (\tilde{w},\tilde{v}) }_{U_p^0}^{\rho-1}),
\end{align*}
$$\norma{F_0(w,v)}_{U^0_q} \leq C (1+ \norma{(w,v)}_{U_p^0}^{\rho}).$$

It remains to ensure that there actually exist $p$ and $q$ that fulfills all the properties \eqref{25.V} to \eqref{25.X}. We already know that, if $X=U_p^0$ and $Y = U^0_q$, Proposition \ref{p22.1} guarantees that, for $\frac{N}{2}<q \leq p$, there exists $\Sigma_{\varphi,0} \subset \rho (-A_0(t))$ and
$$\norma{(\lambda + A_0(t))^{-1}}_{\mathcal{L}(U_q^0,U_p^0)} \leq \dfrac{C}{|\lambda|^{\beta}+1}, \qquad \forall \lambda \in \Sigma_{\varphi,0} \cup \{0\},$$
for any $0< \beta < 1-\frac{N}{2q}-\frac{1}{2} \left( \frac{1}{q} - \frac{1}{p} \right).$ When $q=p$
$$\norma{(\lambda+A_0(t))^{-1}}_{\mathcal{L}(U_p^0)} \leq \dfrac{C}{|\lambda|^{\alpha} +1 },\qquad \forall \lambda \in \Sigma_{\varphi,0} \cup \{0\},$$
for any $0< \alpha < 1-\frac{N}{2p}.$

For $\rho = \frac{p}{q}$ all conditions are satisfied, except \eqref{25.IX}. We only need to verify that there is $\alpha$ and $\beta$ in the previous range such that $1 \leq \rho < \frac{\beta}{1-\alpha}$.

\begin{lemma}\label{l26.1}
	Let $N<q \leq p \mbox{ and }\rho = \frac{p}{q}$.	Then there exist $0<\beta < 1- \frac{N}{2q} - \frac{1}{2} \left( \frac{1}{q}-\frac{1}{p} \right)$ and $0< \alpha < 1- \frac{N}{2p}$ such that $1\leq \rho < \frac{\beta}{1-\alpha}$
	if and only if, for fixed $p>N$, we have $\frac{p(2N+1)}{2p+1}<q \leq p$.

	\begin{proof}
		It is enough to obtain $p$ and $q$ such that  
		$\frac{p}{q}< \frac{1-\frac{N}{2q}-\frac{1}{2}\left(\frac{1}{q} - \frac{1}{p}  \right) }{1- \left( 1-\frac{N}{2p}  \right)}$
		and this is accomplished if, and only if, $ q> \frac{p(2N+1)}{2p+1}$.
%		\begin{align*}
%		\dfrac{1-\frac{N}{2q}-\frac{1}{2}\left(\frac{1}{q} - \frac{1}{p}  \right) }{\frac{N}{2p}} > \dfrac{p}{q} %& \Leftrightarrow q-\dfrac{N}{2} -\dfrac{1}{2} \left( 1- \dfrac{q}{p} \right) > \dfrac{N}{2}\\
%		%& \Leftrightarrow q \left( 1+ \dfrac{1}{2p} \right) > N+\dfrac{1}{2} \\
%		\Leftrightarrow q> \dfrac{p(2N+1)}{2p+1}.
%		\end{align*}
		Note that $q-N> \frac{1}{2}- \frac{q}{2p}>0$, therefore $q>N$ and $\frac{p(2N+1)}{(2p+1)}<p,$ because $\frac{(2N+1)}{(2p+1)}<1$.
	\end{proof}
\end{lemma}

We impose $p$ large enough so that $\alpha$ can be chosen such that
$\alpha > \frac{1}{2} $ and $\alpha+\frac{\delta}{2} > 1.$ Under those hypotheses, the conditions of Theorem \ref{t25.4} are satisfied and we have the existence of local solution.

\begin{proposition}
	
	For $p>N$, $\frac{p(2N+1)}{2p+1} < q \leq p $, $X=U_p^{\varepsilon}$, $Y=U_q^{\varepsilon}$ and the functions $a:\R \rightarrow \R^{+}$, $f:\R \rightarrow \R$ satisfying the conditions mentioned in the beginning of the section, problems \eqref{abs_e} and \eqref{abs_0} have local mild solution.

\end{proposition}

%%%%%%%%%%%%%%%%%%%%%%%%%%%%%%%%%%%%%%

%%%%%%%%%%%%%%%%%%%%%%%%%%%%%%%%%%%%%%

%%%%%%%%%%%%SECTION 4%%%%%%%%%%%%%%%%%%

%%%%%%%%%%%%%%%%%%%%%%%%%%%%%%%%%%%%%%

%%%%%%%%%%%%%%%%%%%%%%%%%%%%%%%%%%%%%%

\section{Positivity and Monotonicity of solutions}
	
	We will prove in this section two features of the parabolic reaction-diffusion equation that allow us to conclude that the solution of the semilinear problem is bounded in bounded intervals. Consequently, this solution is globally defined.

	The first feature is presented in Theorem \ref{t40.1}, which states that a solution that starts in a positive initial condition remains positive as long as it exists. The other is given in Theorem \ref{t40.2}  and it corresponds to the order preserving of the semilinear evolution equation, that is, if $(w_0,v_0) \leq (w_1,v_1)$ in $U_p^{0}$, then the solution obtained for the first initial condition is less than the solution obtained for the second, as long as they exist.
	
	We will follow the ideas presented in \cite[Apendix A]{Arrieta_Carvalho_Anibal} in order to obtained these two theorems mentioned above.
	
	With those properties, we will be able to locate our mild solution 
$(w,v) : (\tau, \tau+\tau_{max}) \rightarrow U_p^{0}$ of \eqref{abs_0} between two bounded functions $\gamma^{+}, \gamma^{-}: I \subset \R \rightarrow U_{p}^{0}$, that is
$$\gamma^{-}(t) \leq (w,v)(t) \leq \gamma^{+}(t), \quad \forall t\in (\tau,\tau+t_{\max}).$$
This implies that the solution $(w,v)(t)$ is bounded.

\begin{definition}
	Given two partially ordered Banach spaces $(X, \leq)$ and $(Y,\leq)$, a map $T:X\rightarrow Y$ is \emph{positive} if $x \geq 0 \Rightarrow Tx \geq 0$ and is \emph{increasing} if $x \leq \tilde{x} \Rightarrow Tx \leq T\tilde{x}$.
\end{definition}

%%%%%%%%%%%%%%%%%%%%%%%%%%%%%%%%%%%%%%

%%%%%%%%%%%%%%%%%%%%%%%%%%%%%%%%%%%%%%

\subsection{The Semigroup of growth $1-\alpha$ is positive}

We already known that, for every $\varepsilon \in [0,1]$ and every 
$t\in \R$, each real number $\lambda >0$ belongs to the resolvent of $-A_{\varepsilon}(t)$. Based in \cite[Theorem 6.43]{Carvalho}, we can prove that this family of operators has the following property: for any $\lambda>0$, $\varepsilon \in[0,1]$ and $t\in \R$
$$U_p^{\varepsilon}  \ni x \geq 0 \Rightarrow (\lambda+A_{\varepsilon}(t))^{-1}x \geq 0$$
(we point out that this include the limiting case and it is a consequence of the fact that the resolvent of $A_{\varepsilon}(t)$ for $\varepsilon \in (0,1]$ converges to the resolvent of $A_0(t)$, see \cite{DDDI}, \cite{DDDIII}).

With the property above, added to the fact that every semigroup can be given by the exponential formula
$$T_{-A_{\varepsilon}(\tau) }(t)x = \lim_{n\rightarrow \infty} \left( I+\dfrac{t}{n} A_{\varepsilon}(\tau) \right)^{-n}x, \quad \forall t>0,$$
it is straight forward to conclude that the semigroup is positive, that is,
$$U_p^{\varepsilon} \ni x \geq 0 \Rightarrow T_{-A_{\varepsilon}(\tau)}(t) x \geq 0 , \quad \forall t>0.$$

%%%%%%%%%%%%%%%%%%%%%%%%%%%%%%%%%%%%%%

%%%%%%%%%%%%%%%%%%%%%%%%%%%%%%%%%%%%%%

\subsection{The Linear process of growth $1-\alpha$ is positive}

We prove now that the linear process $\{  U_{\varepsilon}(t,\tau): U_p^{\varepsilon} \rightarrow U_p^{\varepsilon}; t > \tau   \}$, solution of the nonautonomous linear problem
\begin{equation*}
\begin{cases}
u_t+A_{\varepsilon}(t)u =0, \quad t>\tau, \\
u(\tau) = u_0 \in U_p^{\varepsilon},
\end{cases}
\end{equation*}
is a positive process, for every $\varepsilon \in [0,1]$.

The proof for the sectorial case, where $\varepsilon \in (0,1]$ is very similar to the proof for the almost sectorial case $\varepsilon =0$.  We point out the differences when needed.

If we tried to prove that $\{U_{\varepsilon}(t,\tau) ; t> \tau \}$ is positive by using \eqref{14.IV}, we would not be able to ensure that the integral is positive. To overcome this, we consider the following: given any $\beta >0$, the family $\{ A_{\varepsilon}(t)+ \beta I; t\in \R \}$ is also uniformly (almost) sectorial and generates a semigroup $T_{-(A_{\varepsilon}(\tau)+ \beta I)}(t-\tau)$ solution of
\begin{equation*}
\begin{cases}
u_t+(A_{\varepsilon}(\tau)+\beta I )u=0, \quad t> \tau, \\
u(\tau) = u_0 \in U_p^{\varepsilon}.
\end{cases}
\end{equation*}

Besides, the sectoriality can be chosen to remain the same (same sectors and same constants that bounds the resolvent), once $-\beta I$ just shift the sector of $-A_{\varepsilon}(t)$ to the left.

The family $\{ U_{\varepsilon}(t,\tau) - T_{-(A_{\varepsilon}(\tau)+\beta I)}(t-\tau) ; t>\tau   \}$ is a solution of the problem
\begin{equation*}
\begin{cases}
u_t+A_{\varepsilon}(t)u = -[ A_{\varepsilon}(t)  - A_{\varepsilon}(\tau)  ] T_{-(A_{\varepsilon}(\tau)+ \beta I)}  (t-\tau)  + \beta T_{-(A_{\varepsilon}(t)+\beta I)}(t-\tau), \quad t>\tau\\
u(\tau) = 0
\end{cases}
\end{equation*}
and proceeding as it was done in \cite[pages 25-30]{CDN} (considering adaptations needed), 
we can rewrite the process $\{U_{\varepsilon}(t,\tau); t>\tau\}$ as 
{ \begin{equation*}%\label{new_carc_process}
	\begin{split}U_{\varepsilon}(t,\tau) = & \hspace{0.2cm} T_{-(A_{\varepsilon}(\tau)+\beta I)}(t-\tau)\\
	& +\int_{\tau}^{t} U_{\varepsilon}(t,s)\left\{    [  A_{\varepsilon}(\tau) - A_{\varepsilon} (s)        ]T_{-(A_{\varepsilon}(\tau)+\beta I)}(s-\tau)  + \beta T_{-(A_{\varepsilon}(\tau)+\beta I)}(s-\tau)      \right\} ds.
	\end{split}
	\end{equation*}}

With this new characterization for the process, we can prove it is positive. Before doing so, we enunciate a lemma, consequence of \cite[Remark 3.2]{DDDII}.

\begin{lemma}\label{l38.6}
	There exist constants $M, \tilde{M}, \nu >0$ such that, for $\varepsilon \in (0,1]$
	\begin{equation*}%\label{38.VIII}
	\norma{T_{-A_{\varepsilon} (\tau)} (t) u_0     }_{W^{1,p} (\Omega_{\varepsilon})} \leq M e^{-\nu t} t^{-\frac{1}{2}} \norma{u_0}_{L^{p}(\Omega_{\varepsilon})}, \quad t >0, u_0 \in L^{p}(\Omega_{\varepsilon}) ,
	\end{equation*}
	\begin{equation*}%\label{38.IX}
	\norma{T_{-A_{\varepsilon} (\tau)} (t) u_0     }_{L^{\infty} (\Omega_{\varepsilon})} \leq \tilde{M} e^{-\nu t} t^{-\frac{1}{2}} \norma{u_0}_{L^{p}(\Omega_{\varepsilon})}, \quad t >0, u_0 \in L^{p}(\Omega_{\varepsilon}) ,
	\end{equation*}
	where $M$ and $\tilde{M}$ depends of $N$ and $p$, and $\nu$ comes from the exponential decay of the semigroups generated by $-A_{\varepsilon}(\tau)$. For the case $\varepsilon = 0$, these two inequalities are, for $t>0$ and $u_0 \in L^{p}(\Omega)\times L^{p}(0,1)  $,
	\begin{equation*}
	\norma{T_{-A_{0} (\tau)} (t) u_0     }_{W^{1,p} (\Omega) \times W^{1,p}(0,1)} \leq M e^{-\nu t} t^{-\frac{1}{2}-(1-\alpha)} \norma{u_0}_{L^{p}(\Omega) \times L^{p}(0,1)},
	%{\small \quad t >0, u_0 \in L^{p}(\Omega)\times L^{p}(0,1) }
	\end{equation*}
	\begin{equation*}
	\norma{T_{-A_{0} (\tau)} (t) u_0     }_{L^{\infty} (\Omega)\times L^{\infty}(0,1)} \leq \tilde{M} e^{-\nu t} t^{-\frac{1}{2}-(1-\alpha)} \norma{u_0}_{L^{p}(\Omega) \times L^{p}(0,1)}.
	% \quad t >0, u_0 \in L^{p}(\Omega)\times L^{p}(0,1) 
	\end{equation*}
	
\end{lemma}

\begin{theorem}\label{t39.2}
	Under the previous hypotheses, the linear process $\{  U_{\varepsilon}(t,\tau): U_p^{\varepsilon} \rightarrow U_p^{\varepsilon}; t> \tau   \}$, $\varepsilon \in (0,1]$, is positive. For the limiting case, the process will be positive if we require the additional condition:
	\begin{equation*}%\label{key}
	\alpha + \dfrac{\delta}{2}>1.
	\end{equation*}
	(This hypothesis does not represent any restriction, once we can take $p$ large enough so $\alpha$ is close to 1, as we can see in Proposition \ref{p22.1}).
	
	\begin{proof}
		
		We begin by considering the almost sectorial case $\varepsilon = 0$ and initial condition $u_0 \in D(A^{2}_{0}(\cdot))$ such that $u_0 \geq \gamma I >0$ and $t_0>\tau$.
		%The general case will follow from the density of $D(A^{2}_{0}(\cdot))$ in $U_p^{0}$, the continuity of the linear operator $U_{0}(t,\tau)$ in $U_p^{0}$ and the fact that the positive cone in $U_p^{0}$, that is, the set
		%$$C^{+} = \{ u \in U_{p}^{0}; u \geq 0  \}$$
		%is a closed set.
		
		%		From the continuity of the semigroup restricted to the set $D(A^{2}_{0}(\cdot))$ (see Lemma \ref{l2.6}), we can choose $t_0>\tau$ close enough to $\tau$ so that $T_{-A_{0}(t)}(t_-\tau)u_0 \geq \tilde{\gamma} I> 0$, for all $t\in [\tau,t_0] $. 
		
		Using the same calculus as in \cite[pages 25-30]{CDN}, $U(t) = U_{0} (t_0,t)  $ will be a fixed point of the contraction operator $S:K(t_0,u_0) \rightarrow K(t_0,u_0)$, where 
		\[  K(t_0,u_0) = \left\{  U(\cdot) \in \mathcal{C} \left([  \tau,t_0  ), \mathcal{L}(U_p^{0}    ) \right) ; \sup_{t \in [\tau,t_0)}     (t_0 - t)^{1-\alpha}\norma{U(t)}_{U_p^{0}} <\infty  \right\} \]  and
		\[      
		(SU)(t)   = T_{-A_{0}(t)}(t_0-t)u_0 +\ds\int_{t}^{t_0} U(s) [A_{0}(t)  -A_{0}(s)  ]   T_{-A_{0}(t)}(s-t)u_0 ds.
		\]
		
		By the previous discussion, the same linear process can be given as a fixed point, in $K(t_0,u_0)$, of the map
		{	\[
			\begin{split}(S^{\beta}U)(t) = & T_{-(A_{0}(t)+\beta I)}(t_0-t)u_0\\
			&+\int_{t}^{t_0} U(s)\left\{    [  A_{0}(t) - A_{0} (s)        ]T_{-(A_{0}(t)+\beta I)}(s-t) u_0 + \beta T_{-(A_{0}(t)+\beta I)}(s-t)u_0      \right\} ds,\end{split}
			\]		}for any $\beta >0.$
		
		By choosing $t_0$ close enough to $\tau$, we have:
		\begin{enumerate}
			\item $ [  A_{0}(t) - A_{0} (s)        ]T_{-(A_{0}(t)+\beta I)}(s-t) u_0$ is bounded in $L^{\infty}(\Omega)\oplus L^{\infty}(0,1)$ and this bound is uniform for $\tau \leq t \leq s \leq t_0$.
			\item There exists $\tilde{\gamma}>0$ such that $T_{-(A_{0}(t)+\beta I)}(s-t)u_0 \geq \tilde{\gamma}I >0$ for all $\tau \leq t \leq s \leq t_0$.
		\end{enumerate}
		
		%	For $u_0 \in D(A^{2}_{0}(\cdot))$, it follows from Lemma \ref{l2.6} that
		%	\[   
		%	[A_{0}(t)- A_{0}(s)] T_{  -( A_{0}(t) + \beta I  )  } (s-t) u_0 = T_{  -( A_{0}(t) + \beta I  )  } (s-t)  [ A_{0} (t) - A(0)(s)   ]u_0
		%	\]  
		%	and is a consequence of Lemma \ref{l38.6} that this is bounded for $s\in [\tau,\tau+t_0]$ in the $L^{\infty}(\Omega_{\varepsilon})$. Once $T_{ -( A_{\varepsilon} (\tau)  +\beta I  ) } (s-\tau)u_0$ is positive and continuous, and choosing a small $t_0$ if necessary, we can garantee that 
		%	\[
		%	T_{ -( A_{\varepsilon} (\tau)  +\beta I  ) } (s-\tau)u_0 \geq \tilde{\gamma} I >0.
		%	\]
		%	
		We will prove those two claims in the sequence, but assuming they are true, we can choose $\beta$ large enough that  
		\[   [  A_{0}(t) - A_{0} (s)        ]T_{-(A_{0}(t)+\beta I)}(s-t)u_0  + \beta T_{-(A_{0}(t)+\beta I)}(s-t) u_0   \geq 0     , \quad \forall \tau \leq t \leq s \leq t_0,\]
		and, by defining 
		\[
		K^{+}(t_0,u_0) = \{   U \in K(t_0,u_0); \mbox{ } U(t) \mbox{ is positive}\},
		\]
		we have $ S^{\beta}|_{K^{+}(t_0,u_0)}:K^{+}(t_0,u_0) \rightarrow K^{+}(t_0,u_0)$. 
		%We assure that the fixed point of $S^{\beta} $ belongs to $K^{+}(t_0,u_0)$. 
		Indeed, for $t = \tau$, $U_{0}(\tau,\tau) = I$ (positive) and once
		\[ T_{-(A_{0}(t)+\beta I)}(t_0-t)u_0 \geq 0,\]
		and
		\[ \int_{t_0}^{t} U(s)    \left\{     [ A_{0}(t)-A_{0}(s) ] T_{-(A_{0}(t)+\beta I)}(s-t) u_0    + \beta T_{  -( A_{0}(t)    +\beta I  )    }(s-t)u_0        \right\} ds   \]
		is continuous, it follows that for small $t_0$, $S^{\beta}|_{K^{+}(t_0,u_0)}:K^{+}(t_0,u_0) \rightarrow K^{+}(t_0,u_0)$. From the uniqueness of the fixed point, we then have that it belongs to $K^{+}(t_0,u_0)$ and then $ U_{0}(t_0,t) $ is positive.
		
		Up to now we can conclude that 
		$$u_0 \geq \gamma I >0 \mbox{ and } u_0 \in D(A^{2}_{0}(\cdot)) \Rightarrow U_{0}(t,\tau)u_0 \geq 0 , \quad \forall t>\tau.$$
		
		If $u_0 \in D(A^{2}_{0}(\cdot))$ with $u_0 \geq 0$, we consider $u_n = u_0+\frac{1}{n}I$. For each $n$, $u_n$ satisfies the previous conditions and $U_{0}(t,\tau)u_n \geq 0$, for every $t>\tau$. From the fact that $u_n \rightarrow u_0$ in $U_p^{0}$ and from the continuity of the linear operator $U_{0}(t,\tau)$, it follows that 
		\[   U_{0}(t,\tau)u_n \rightarrow  U_{0}(t,\tau)u_0 \mbox{ in } U_p^{0}   \Rightarrow U_{0}(t,\tau)u_0 \geq 0.  \]
		
		For the general case where $u_0 \in U_p^{0}$ with $u_0 \geq 0$ we use the density of $D(A^{2}_{0}(\cdot))$ to obtain $(u_n)\subset D(A^{2}_{0}(\cdot))$, $u_n \geq 0$, such that $u_n \rightarrow u_0$ in $U_p^{0}$. Then, $U_{0}(t,\tau) u_0 \geq 0$.

		Therefore, the process $\{  U_{0}(t,\tau):U_p^{0} \rightarrow U_p^{0}; t>\tau   \}$ will be positive. We only need to prove the two previous claims.
		\medskip
		
		\textbf{Claim 1:} %For $t_0-\tau>0$ small enough and $u_0 \in D(A^{2}_0(\cdot))$, $ [  A_{0}(t) - A_{0} (s)        ]T_{-(A_{0}(t)+\beta I)}(s-t) u_0$ is bounded in $U_{\infty}^{0}=L^{\infty}(\Omega)\oplus L^{\infty}(0,1)$ and this bound being uniform for $\tau \leq t \leq s \leq t_0$.
		Fixing any $t^{*} \in \R$, $u_0 \in D(A_0^{2}(t^{*}))$ we have $ u_0 \in D(A^{2}_0(t^{*})+ \beta I)$. Therefore, if 
		$$ w_0 = (A_0(t^{*})+\beta I)u_0, $$
		then $w_0 \in D(   A_0(t^{*})    + \beta I     )$ and $u_0 = \{   A_0(t^{*})    +\beta I   \}^{-1} w_0$. So
		\begin{equation}\label{step1}
		\begin{split}
		& 	[  A_{0}(t) - A_{0} (s)        ]T_{-(A_{0}(t)+\beta I)}(s-t) u_0 \\
		& = [  A_{0}(t) - A_{0} (s)        ]   \{  A_0(t)    +\beta I      \}^{-1} \{  A_0(t)    +\beta I      \} T_{-(A_{0}(t)+\beta I)}(s-t) \{   A_0(t^{*})    +\beta I   \}^{-1} w_0\\
		& = [  A_{0}(t) - A_{0} (s)        ]   \{  A_0(t)    +\beta I      \}^{-1}  T_{-(A_{0}(t)+\beta I)}(s-t)\{  A_0(t)    +\beta I      \} \{   A_0(t^{*})    +\beta I   \}^{-1} w_0\
		\end{split}
		\end{equation}
		due to Lemma \ref{l2.6} adapted to the family $\{A_0(\cdot)+ \beta I\}$. If $\xi$ is an element of $D(A_0(\cdot))$, 
		\[ \norma{  [A_0(t)    -A_0(s)]    (A_0(t)+ \beta I)^{-1}  \xi  }_{U_{\infty}^{0}} \leq C |t-s|^{\delta}  \norma{\xi}_{U_{\infty}^{0}} .    \]
		
		Indeed, it is straightforward that $\{ A_0(\cdot)  + \beta I  \}$ is H\"{o}lder continuous and for $p>N$ 
		\[    \norma{  [A_0(t)    -A_0(s)]    (A_0(t)+ \beta I)^{-1}  \xi  }_{U_{p}^{0}} \leq C |t-s|^{\delta}  \norma{\xi}_{U_{p}^{0}}.   \]

		From the fact that $ \xi \in D(A_0(\cdot))   \hookrightarrow W^{{2,p}} (\Omega) \oplus W^{2,p}(0,1)  \hookrightarrow L^{\infty}(\Omega) \oplus L^{\infty}(0,1),  $ it follows that 
		\[   \norma{\xi}_{U_{p}^{0}} \leq K \norma{\xi}_{U_{\infty}^{0}} \Rightarrow  \norma{  [A_0(t)    -A_0(s)]    (A_0(t)+ \beta I)^{-1}  \xi  }_{U_{p}^{0}} \leq C |t-s|^{\delta}  \norma{\xi}_{U_{\infty}^{0}}.   \]  
		
		Now, we have $[A_0(t)    -A_0(s)]    (A_0(t)+ \beta I)^{-1}  \xi \in D(A_0(\cdot)) \hookrightarrow L^{\infty}(\Omega) \oplus L^{\infty} (0,1)$ and
		\begin{align*}
		&	\norma{  [A_0(t)    -A_0(s)]    (A_0(t)+ \beta I)^{-1}  \xi  }_{U_{\infty}^{0}} \\
		& = \limsup_{p \rightarrow \infty} \norma{  [A_0(t)    -A_0(s)]    (A_0(t)+ \beta I)^{-1}  \xi  }_{U_{p}^{0}} \leq C|t-s|^{\delta} \norma{\xi}_{U_{\infty}^{0}} .
		\end{align*}

		Returning to \eqref{step1} we have
		\begin{equation}\label{step2}
		\begin{split}
		& \norma{   	[  A_{0}(t) - A_{0} (s)        ]   \{  A_0(t)    +\beta I      \}^{-1}  T_{-(A_{0}(t)+\beta I)}(s-t)\{  A_0(t)    +\beta I      \} \{   A_0(t^{*})    +\beta I   \}^{-1} w_0   }_{    U_{\infty}^{0}   } \\
		& \leq C |t-s|^{\delta}    \norma{  T_{-(A_{0}(t)+\beta I)}(s-t)\{  A_0(t)    +\beta I      \} \{   A_0(t^{*})    +\beta I   \}^{-1} w_0     }_{U_{\infty}^{0}}.
		\end{split}
		\end{equation}
		
		Note that Lemma \ref{l38.6} remains valid for the semigroup $T_{-( A_0(t) + \beta I   )}$ and, as a consequence, 
		\begin{align*}
		& \norma{  T_{-(A_{0}(t)+\beta I)}(s-t)\{  A_0(t)    +\beta I      \} \{   A_0(t^{*})    +\beta I   \}^{-1} w_0     }_{U_{\infty}^{0}}\\
		& \leq M e^{-\nu (s-t)}  (s-t)^{-\frac{1}{2} - (1-\alpha)} \norma{ \{  A_0(t)    +\beta I      \} \{   A_0(t^{*})    +\beta I   \}^{-1} w_0    }_{U_p^{0}}\\
		& \leq M e^{-\nu (s-t)}  (s-t)^{-\frac{1}{2} - (1-\alpha)} \norma{w_0    }_{U_p^{0}},
		\end{align*}
		since $ \{  A_0(t)    +\beta I      \} \{   A_0(t^{*})    +\beta I   \}^{-1} $ is bounded in $\mathcal{L}(U_p^{0})$. Therefore, $$ T_{-(A_{0}(t)+\beta I)}(s-t)\{  A_0(t)    +\beta I      \} \{   A_0(t^{*})    +\beta I   \}^{-1} w_0  $$
		is bounded in $L^{\infty}(\Omega) \oplus  L^{\infty}(0,1)$ if $s$ is not near the time $t$. For the case it is, we will use the facts that $u_0 \in D(A^{2}_0(\cdot))$, $w_0 \in D(A_0(\cdot))$ and the continuity of the semigroup for these cases.
		
		We know that 
		\[ \norma{\{ A_0(t) +\beta I  \}      \{  A_0(t^{*})+\beta I \}^{-1}    w_0   }_{L^{\infty}(\Omega)  \oplus  L^{\infty}(0,1)}   \leq C \norma{w_0}_{L^{\infty}(\Omega)  \oplus  L^{\infty}(0,1)}  \]
		and if $\xi(t) = \{ A_0(t) +\beta I  \}      \{  A_0(t^{*})+\beta I \}^{-1}    w_0   $, the previous inequality allows us to conclude that
		\[    \xi(t) \in D(A_0(\cdot))  \quad \mbox{and} \quad  \xi(t) \in B_{U_{\infty}^{0}} \left[0, C \norma{w_0}_{U_{\infty}^{0}}\right], \forall t \in [\tau, t_0).  \]
		
		From Lemma \ref{l2.6}, it follows that $$   \norma{ T_{ -(A_0(t) + \beta I) }  (s-t) \xi(t)  }_{U_p^{0}} \stackrel{s\rightarrow t^{+}}{\longrightarrow}  \norma{ \xi(t) }_{U_p^{0}}.  $$
		
		On the other hand, $ T_{ -(A_0(t) + \beta I) }  (s-t) \xi(t)   \in D(A_0(\cdot))$ for every $s \geq t$. Therefore, belongs to $U_{\infty}^{p}.$ Making $p \rightarrow \infty$ and using an argument of choosing a diagonal sequence, we have
		$$   \norma{ T_{ -(A_0(t) + \beta I) }  (s-t) \xi(t)  }_{U_{\infty}^{0}} \stackrel{s\rightarrow t^{+}}{\longrightarrow}  \norma{ \xi(t) }_{U_{\infty}^{0}} $$
		for every $t \in [\tau, t_0)$. Therefore, 
		\[  \norma{ T_{-(A_{0}(t)+\beta I)}(s-t)\{  A_0(t)    +\beta I      \} \{   A_0(t^{*})    +\beta I   \}^{-1} w_0  }_{U_{\infty}^{0}}  \in B_{U_{\infty}^{0}} \left[0,2C \norma{w_0}_{U_{\infty}^{0}}\right]     \]
		for $s$ close to $t$ and any $t \in [\tau, t_0)$.
		
		It follows from \eqref{step2} that the Claim 1 is valid.

		\mbox{ }
		
		\textbf{Claim 2:} We prove that there exists $t_0>\tau$ such that $T_{-A_0(t)}(t-\tau)u_0 \geq \tilde{\gamma} I >0$ for all $\tau \leq t <t_0$. From Lemma \ref{l2.6} we already know that
		$$T_{-A_0(\tau)}(t-\tau)u_0 \stackrel{t \rightarrow \tau^{+}}{\longrightarrow} u_0$$
		and then, using similar arguments, we can choose $t_0-\tau>0$ small so that $T_{-A_0(\tau)}(t-\tau)u_0 \geq \tilde{\gamma}I>0$. It follows from \cite[Lemma 2.2]{CDN} that, for $\alpha+\frac{\delta}{2}-1 >0$,
		$$ \norma{ T_{-A_0(t)} (t-\tau)  u_0    - T_{-A_0(\tau)}   (t-\tau)  u_0   }_{U_p^{0}}  \leq C (t-\tau)^{-2+2\alpha + \delta}\norma{u_0}_{U_p^{0}}. $$
		From this and the same reasoning used before (letting $p \rightarrow \infty$), we conclude that the two norms $ \norma{ T_{-A_0(t)} (t-\tau)  u_0 }_{U_p^{\infty}} $ and $ \norma{T_{-A_0(\tau)}   (t-\tau)  u_0   }_{U_p^{\infty}}$ are close. This ensure the existence of $t_0>\tau$ such that Claim 2 holds true.

		\mbox{ }
		
		The case $\varepsilon \in (0,1]$ follows using the same idea, but in order to prove that 
		$$T_{-(A_{\varepsilon}(t)+\beta I)}(s-t)u_0 \geq \tilde{\gamma}I >0$$ for all $\tau \leq t \leq s \leq t_0$,
		we only need $\alpha + \delta>1$. Besides that, the chosen space $K(t_0,u_0)$ will be $ K(t_0,u_0) =\{   U(\cdot) \in \mathcal{C}(\tau,t_0]; \sup_{ t\in [\tau, +t_0) } \norma{U(s)   } <\infty  \}.$
		
	\end{proof}
\end{theorem}

%%%%%%%%%%%%%%%%%%%%%%%%%%%%%%%%%%%%%%

%%%%%%%%%%%%%%%%%%%%%%%%%%%%%%%%%%%%%%

\subsection{The solution of the semilinear equation is positive}

Consider the family $\{  A_{\varepsilon}(t)   +\beta I ; t\in \R \} $ of uniformly (almost) sectorial operators and  denote by $\{ U_{\varepsilon}(t,\tau):U_p^{\varepsilon} \rightarrow U_p^{\varepsilon}; t > \tau  \}$ and $\{ U_{\varepsilon}^{\beta}(t,\tau):U_p^{\varepsilon} \rightarrow U_p^{\varepsilon}; t >\tau  \}$ the family of operators solution of
$$\begin{matrix}
\begin{cases}
u_t +A_{\varepsilon}(t)u =0,\\
u(\tau) = u_0 \in U_p^{\varepsilon},
\end{cases} \quad \quad & \quad \quad \begin{cases}
u_t +[A_{\varepsilon}(t)+\beta I]u =0,\\
u(\tau) =u_0 \in U_p^{\varepsilon},
\end{cases} 
\end{matrix}
$$
respectively. With same reasoning used in the previous section, the process $  \{ U_{\varepsilon}^{\beta}(t,\tau);t>\tau  \}$ is positive.

We know that the semilinear evolution equations \eqref{abs_e} and \eqref{abs_0}
%$$	\begin{cases}
%u_t+A_{\varepsilon}(t) u = F_{\varepsilon}(u), \quad t> \tau, \\
%u(\tau) = u_0 \in U_p^{\varepsilon}
%\end{cases}$$
have a unique local mild solution which is the fixed point of 
$$  \mathcal{F_{\varepsilon}}(u)(t)   = U_{\varepsilon}(t,\tau)u_0 + \int_{\tau}^{t}  U_{\varepsilon} (t,s)   F_{\varepsilon}(u(s))  ds, \quad \varepsilon \in [0,1].$$

The operator $\mathcal{F_{\varepsilon}}$ is a contraction on the Banach space $K_{\varepsilon}(t_0,u_0)$, 
$$K_{\varepsilon}(t_0,u_0) = \left\{ u \in \mathcal{C}( (\tau, t_0], U_p^{\varepsilon}  ) ;   \norma{u}_{L^{\infty} ([\tau, t_0], U_p^{\varepsilon})   }   \leq \mu +C \norma{u_0}_{U_p^{\varepsilon}}  \right\}, \quad  \mbox{for }\varepsilon \in (0,1],$$
$$ K_{0}(t_0,u_0) =  \left\{ u \in \mathcal{C}( (\tau, t_0], U_p^{0}  ) ;   \norma{u(t)}_{ U_p^{0}  }   \leq (t-\tau)^{\alpha-1} \left( \mu +C \norma{u_0}_{U_p^{0}}\right)  \right\}, \quad \mbox{for }\varepsilon = 0,$$
and $\mu>0$ is chosen properly to assure $\mathcal{F_{\varepsilon}}$ is a contraction (see \cite[pages 34,35]{CDN}).

We can now prove the positivity of solutions of the semilinear equations.

\begin{theorem}\label{t40.1}
	For $\varepsilon \in (0,1]$, if $u_0 \in U_p^{\varepsilon}$ is such that $u_0 \geq 0$, then the solution $u(t,\tau,u_0)$ of the semilinear equation
	$$\begin{cases}
	u_t+A_{\varepsilon}(t) u = F_{\varepsilon}(u), \quad t> \tau, \\
	u(\tau) =u_0 \in U_p^{\varepsilon},
	\end{cases}
	$$ is positive as long as it exists. 
	
	For the limiting case $\varepsilon =0$, assuming that $\alpha + \frac{\delta}{2}>1$, for every $u_0 \in D(A_0(\cdot))$ such that $u_0 \geq 0$, the solution of 
	$$\begin{cases}
	u_t+A_{0}(t) u = F_0(u), \quad t>\tau,\\
	u(\tau) =u_0 \in D(A_0(\cdot)),
	\end{cases}
	$$  is positive as long as it exists.

\begin{proof}
	We deal with the case $\varepsilon=0$, since the nonlimiting case is simpler. We start by considering $u_0 \in D(A_0(\cdot))$ such that $u_0 \geq \gamma  >0$. We know that the mapping 
	$$ \mathcal{F}_0(u)(t) = U(t,\tau)u_0 + \int_{\tau}^{t} U(t,s) F(u(s))ds$$
	is a contraction in $K_0(t_0,u_0)$.
	Once $u_0 \in D(A_0(\cdot))$, for $t_0-\tau>$ small, the function $u(\cdot) \in \mathcal{C}( [\tau, t_0); U_{\infty}^{0}  )$ and we can assume $u(t) \geq \tilde{\gamma}  >0$ for some $\tilde{\gamma}>0$. Therefore, $F(u(s))$ is bounded for $s\in [\tau,t_0).$
	
	Let $\beta >0$ be such that 
$F(u(s))+ \beta  u(s) \geq 0,$ for $s \in [\tau,t_0),$
	and consider the problem
	$$\begin{cases}
	u_t+[A_0(t)+\beta I ]u = F(u)+\beta  u, \quad t>\tau,\\
	u(\tau) =u_0 \in D(A_0(\cdot)).
	\end{cases}$$
	
	Clearly, $u(t,\tau,u_0)$ is also a solution  for this problem, but now it can be obtained as a fixed point of
	$$  \mathcal{F}_{0}^{\beta}(u)(t)   = U_{0}^{\beta}(t,\tau)u_0 + \int_{\tau}^{t}  U^{\beta}_{0} (t,s)[F_{0}(u(s)) + \beta  u(s)] ds.$$
	
	From the previous discussion, $\{U_0^{\beta}(t,\tau); t> \tau\}$ is also positive and $F(u(s)) + \beta  u(s) \geq 0$. Therefore, if we restrict the map $\mathcal{F}_0^{\beta}$ to 
	$$K_0^{+}(t_0,u_0) = \{ u \in K_0(t_0,u_0) ; \mbox{ }u \geq 0 \},$$
	we have $\mathcal{F}_{0}^{\beta}(K_0^{+}(t_0,u_0))\subset K_0^{+}(t_0,u_0)$. Hence, the fixed point of $\mathcal{F}_0^{\beta}$ belongs to $K_0^{+}(t_0,u_0)$, that is, $u(t,\tau,u_0) \geq 0$, for all $ t \in [\tau, t_0)$.
	
	For the general case where $u_0 \in D(A_0(\cdot))$ with $u_0 \geq 0$, we consider the sequence $u_n = u_0 + \frac{1}{n}I$ and let $u_n(t)$ be the solution in its maximal interval $[\tau, t_n)$ of the problem
	$$
	\begin{cases}
	u_t + A_0(t)u = F_0(u), \quad t>\tau,\\
	u(\tau) = u_n.
	\end{cases}
	$$
	
	Once we have $u_n \rightarrow u_0$ in $U_p^{0}$ and in $U_{\infty}^{0}$, and from the fact that the maximal interval of existence can be taken uniformly in bounded set (see discussion after Theorem \ref{t25.4}), we have $t_n =t_0$. It follows from Proposition \ref{t25.5}, item (2), that $u_n(t) \rightarrow u_0(t)$ in $U_p^{0}$, for $t \in [\tau,t_0)$, which guarantees that $u_0(t) \geq 0$.
	
	With continuation arguments, we conclude that the solution is positive as long as it exists.

	For the nonlimiting case $(\varepsilon \in (0,1])$, we have the continuity of the solution for every initial condition $u_0 \in U_p^{\varepsilon}$. Therefore we can ensure that $F(u(s))$ is bounded for every $s\in [\tau,t_0)$. Besides, we do not need the hypothesis $\alpha+\frac{\delta}{2}>1$ to ensure the process $ \{U_{\varepsilon}^{\beta}(t,\tau); t> \tau\}  $ is positive.
\end{proof}
\end{theorem}

\begin{remark}\label{r1}
	Even though positive initial conditions $u_0 \in U_p^{0}$ imply that 
	$$T_{-A_0(\tau)}(t)u_0 \geq 0 \quad \mbox{and} \quad U_0(t,\tau)u_0 \geq 0,$$
	this, in general, will not be true for the solution $u(t,\tau,u_0)$ of the semilinear equation. This is a consequence of the discontinuity that might appear in the initial time, added to the presence of the nonlinearity $F_0$.
	
	However, Theorem \ref{t40.1} assures the positivity of $u(t,\tau,u_0)$ when the initial condition is in the domain of $A_0(\tau)$. For this case, we remove the discontinuity on the initial time.
	
	This is an important difference between the limiting case and the case $\varepsilon \in (0,1]$, where the positivity is assured for every initial condition $u_0 \in U_p^{\varepsilon}$.
\end{remark}

%%%%%%%%%%%%%%%%%%%%%%%%%%%%%%%%%%%%%%

%%%%%%%%%%%%%%%%%%%%%%%%%%%%%%%%%%%%%%

\subsection{Comparison results}

We now prove the order preserving property of the semilinear evolution equation. This is done in the next theorem and we compare solutions under two different scenarios: One when we make changes in the initial condition and other when we make changes in the nonlinearity.

\begin{theorem}\label{t40.2}
	Let $F_{\varepsilon}(\cdot), G_{\varepsilon}(\cdot): U_p^{\varepsilon} \rightarrow U_p^{\varepsilon}$ be maps satisfying the following hypothesis:
	for every $r>0$, there exists $\beta = \beta (r)$ such that $F_{\varepsilon}(\cdot) + \beta I$ and $G_{\varepsilon}(\cdot) + \beta I$ are increasing in $B_{U_p^{\varepsilon}}(0,r)$. For $\varepsilon \in(0,1]$, we have
	\begin{enumerate}
		\item If $u_0, u_1 \in U_{p}^{\varepsilon}$, $u_0 \geq u_1$, then the solutions $u(t,\tau,u_i)$, $i=0,1$, of the problems 
		\begin{equation*}
		u_t+A_{\varepsilon}(t) u = F_{\varepsilon}(u); \quad
		u(\tau)=u_i,
		\end{equation*}
		will satisfy $u(t,\tau,u_0 ) \geq u(t,\tau,u_1)$ as long as they exist.
		\item If $F_{\varepsilon}(\cdot)   \geq   G_{\varepsilon}(\cdot)$ for every $u_0\in U_p^{{\varepsilon}}$, and $u_F(t,\tau, u_0)$, $u_G(t,\tau,u_0)$ are the solutions of the evolution equations 
		%with nonlinearity $F_{\varepsilon}$ and $G_{\varepsilon}$, 
		\begin{align*}
		u_t+A_{\varepsilon}(t) u & = F_{\varepsilon}(u); \quad
		u(\tau)=u_0, \\
		u_t+A_{\varepsilon}(t) u & = G_{\varepsilon}(u); \quad
		u(\tau)=u_0,
		\end{align*}
		respectively, we have $u_F(t,\tau,u_0) \geq u_G(t,\tau,u_0)$ as long as they exist.
	\end{enumerate}
	
	For the limiting case $\varepsilon =0$, the two claims above remain valid if $u_0, u_1 \in D(A_0(\cdot))$ and $\alpha +\frac{\delta}{2}>1$.

\begin{proof}
	We prove the theorem for the limiting case. Let $u_0,u_1 \in D(A_0(\cdot))$. We already know that $u_i(t) = u_F(t,\tau,u_i)$ is the unique fixed point of 
	$$ \mathcal{F}_0^{\beta}(u)(t) = U_0^{\beta} (t,\tau)u_i + \int_{\tau}^{t} U_0^{\beta} (t,s) \left[  F_0(u(s))  + \beta I u(s)    \right] ds$$
	in the space $$V_i =\left\{  u \in \mathcal{C}( (\tau,t_0) , U_p^{0}  ); \quad u(\tau) = u_i; \quad \norma{u}_{  L^{\infty} (( \tau, t_0)    , U_p^{0} ) } \leq (t-\tau)^{\alpha-1} (\mu + C\norma{u_i}_{U_p^{0}}   )        \right\},$$ for small $t_0-\tau >0$.
	
	Consider $ V_0^{+} = \{ u \in V_0; u(t) \geq u_1(t)    \}.$
	This space is not empty for small $t_0-\tau$, once $u_1(t) + u_0 -u_1$ belongs to $V_0^{+}$ (due to the continuity of the functions $u_0(t)$ and $u_1(t)$ at $[\tau,t_0)$).
	%, since the initial conditions belong to the domain of $A_0(\cdot)$).
	
	Besides, $F_0(u_0)+ \beta I u_0 \geq F_0(u_1) + \beta I u_1$ for $\beta$ large enough and this map is increasing. Hence, $ \mathcal{F}_0^{\beta}(V_0^{+}) \subset V_0^{+}.$
	Therefore, the fixed point of $\mathcal{F}_0^{\beta}$ belongs to $V_0^{+}$ and $u_0(t) \geq u_1(t),$ for $t \in [\tau, t_0)$.
	
	A continuation argument implies that $u_0(t) \geq u_1(t),$ as long the solutions exist.

	The case where we change the nonlinearities follows the same ideas.
	
	As for the case $\varepsilon \in (0,1]$, we proceed the same way, without the necessity of imposing that the initial conditions to belong in the domain of $A_{\varepsilon}(\cdot)$.%, once we have continuity of solutions in the initial time. 
\end{proof}
\end{theorem}

\begin{remark}
	As a consequence of Theorems \ref{t40.1}, \ref{t40.2}  and from what we have pointed out in Remark \ref{r1}, we will only have comparison results for the limiting case when the initial conditions belong to $D(A_0(\cdot))$.
\end{remark}

%%%%%%%%%%%%%%%%%%%%%%%%%%%%%%%%%%%%%%

%%%%%%%%%%%%%%%%%%%%%%%%%%%%%%%%%%%%%%

%%%%%%%%%%%% SECTION 5 %%%%%%%%%%%%%%%%%%

%%%%%%%%%%%%%%%%%%%%%%%%%%%%%%%%%%%%%%

%%%%%%%%%%%%%%%%%%%%%%%%%%%%%%%%%%%%%%

\section{Global solutions and existence of attractors}

The previous results about positivity and comparison of solutions are  applied in this section to prove Theorem \ref{t41.1,41.2,41.4} below, which ensures that the solutions obtained for the problem are globally defined. To this end, we impose that the nonlinearity $f: \R \rightarrow \R$ satisfies the dissipativeness condition
\begin{equation}\label{D}
\limsup_{|s| \rightarrow \infty} \dfrac{f(s)}{s} < 1.
\end{equation}
%that is, this limsup will be less than the first eigenvalue of $ A_0(t)$. For further details on the spectrum of this operator, we recommend \cite[Section 3.2]{DDDII}. We just mention here that, for each $t \in \R$, the spectrum is given by
%$$\sigma (A_0(t)) = \{ \lambda_i ^{t}\in \R; i\in \N, \lambda_i^{t} \mbox{ is an eigenvalue and } 1 \leq \lambda_1^{t} \leq \lambda_2^{t} \leq ... \leq \lambda_n^{t} \leq ...   \}.$$

This dissipativeness implies the next result.

\begin{lemma}\label{dis_cond}
	There exist $r>0$ and $\xi \in (0,1)$ such that, for every $s>r$, $f(s)<\xi s.$
	
	\begin{proof}
		Let $\xi := \sup_{|s|>r} \frac{f(s)}{s} <1 $. 
		%$ \inf_{r>0} \left\{   \sup_{|s|> r}   \frac{f(s)}{s}     \right\} <1$. So, 
		There exists $r>0$ such that, for all $|s|>r$, $\frac{f(s)}{s} \leq \xi$. Therefore, if $s> r$, than $f(s)\leq \xi s $.
%		$$\xi := \sup_{|s|>r} \frac{f(s)}{s} <1 \Rightarrow \forall |s|>r, \quad \dfrac{f(s)}{s} \leq \xi \Rightarrow \forall s > r , \quad f(s)\leq \xi s .$$
		%		If $u \in U_{\infty}^{\varepsilon}$ is a positive element with $u(\cdot) \geq r_0$ in $\Omega_{\varepsilon}$, then
		%		\begin{equation}
		%		F_{\varepsilon}(u)\leq \xi_0 u
		%		\end{equation}
	\end{proof}
\end{lemma}

For $s \in [0,r]$, there exists $M>0$ such that $f(s) \leq M$, since $f$ is continuous. Then
\begin{equation*}
f(s) \leq \begin{cases}
\xi s, \mbox{ if } s>r, \\
M, \mbox{ if }s\in [0,r],
\end{cases}
\end{equation*}
which implies that, if $v$ is a positive solution in $U_{p}^{\varepsilon}$,
\begin{equation*}
F_{\varepsilon}(v)(x) = f(v(x))  \leq \begin{cases}
\xi v(x), \mbox{ if } v(x)>r ,\\
M, \mbox{ if }v(x)\in [0,r].
\end{cases}
\end{equation*}
By considering in $U_p^{\varepsilon}$ the constant function $K_M(x) = M$, we have
\begin{equation*}%\label{comparison_no_linearities}
F_{\varepsilon}(\cdot) \leq \xi \cdot + K_M.
\end{equation*}

To prove the mild solutions are globally defined in time we focus on the case $\varepsilon =0$.
Consider, for $u_0 \in D(A_0(\cdot))$, the three following evolution equations, for $t> \tau$,
$$\begin{matrix}
\begin{cases}
u_t+A_0(t)u=F_{0}(u), \\
u(\tau) = u_0 \in D(A_0(\cdot)),
\end{cases}
&
\begin{cases}
u^{+}_t+A_0(t)u^{+}=F_{0}u^{+} ,\\
u^{+}(\tau) = |u_0| \in D(A_0(\cdot)),
\end{cases}
&
\begin{cases}
\gamma_t+A_0(t)\gamma=\xi \gamma + K_M,\\
\gamma (\tau) = |u_0|  \in D(A_0(\cdot)),
\end{cases}
\end{matrix}$$
with solutions, respectively, given by
$$u(\cdot, \tau, u_0), \mbox{ } u^{+}(\cdot,\tau, |u_0|) \mbox{ and } \gamma (\cdot, \tau, |u_0|),$$
defined each on its maximal interval. Then we use the comparison results proved previously to conclude that 
\begin{equation*}
u(t, \tau, u_0)\leq u^{+}(t,\tau, |u_0|) \leq  \gamma (t, \tau, |u_0|),
\end{equation*}
as long as they exist. In the same way, if we consider the problems, for $t> \tau$,
\[
\begin{cases}
u_t+A_0(t)u=F_{0}(u), \\
u(\tau) = u_0 \in D(A_0(\cdot)),
\end{cases}
\begin{cases}
u^{-}_t+A_0(t)u^{-}=F_{0}u^{-} ,\\
u^{-}(\tau) = -|u_0| \in D(A_0(\cdot)),
\end{cases}
\begin{cases}
\gamma^{-}_t+A_0(t)\gamma^{-}=\xi \gamma^{-} -K_M,\\
\gamma^{-} (\tau) =- |u_0|  \in D(A_0(\cdot)),
\end{cases}
\]
we would have
\begin{equation*}
u(t, \tau, u_0)\geq u^{-}(t,\tau, -|u_0|) \geq  \gamma^{-} (t, \tau, -|u_0|),
\end{equation*}
where these functions are the solutions of the above problems, respectively, and the inequality is valid as long as they exists. Due to the linearity and the uniqueness of solution, we can easily show that $\gamma^{-}(t,\tau,-|u_0|) = - \gamma (t,\tau,|u_0|)$. Therefore,
\begin{equation}\label{40.IV}
|u(t,\tau, u_0) | \leq \gamma (t,\tau,|u_0|),
\end{equation}
as long as the solutions exist. We have located the mild solution for the semilinear equation among two functions $\gamma$ and $-\gamma$ which are now solutions of linear equations and, therefore, is an easy task to obtain bounds for them. To do it, we need the following lemma, which is a consequence of Lemma \ref{l38.6}, the bounds \eqref{31.IV}, \eqref{31.V} and the formulation \eqref{31.III} that we have for the process.

\begin{lemma}\label{l40.3}
	There exist $M, \tilde{M}, \nu >0$ such that, for $\varepsilon \in (0,1]$
	\begin{equation*}%\label{40.V}
	\norma{U_{\varepsilon} (t,\tau)u_0 }_{W^{1,p} (\Omega_{\varepsilon})} \leq M e^{-\nu (t-\tau)}(t-\tau)^{-\frac{1}{2}} \norma{u_0}_{L^{p}(\Omega_{\varepsilon})}, \quad t >\tau, u_0 \in L^{p}(\Omega_{\varepsilon}), 
	\end{equation*}
	\begin{equation*}%\label{40.VI}
	\norma{U_{\varepsilon} (t,\tau) u_0     }_{L^{\infty} (\Omega_{\varepsilon})} \leq \tilde{M} e^{-\nu (t-\tau)} (t-\tau)^{-\frac{1}{2}} \norma{u_0}_{L^{p}(\Omega_{\varepsilon})}, \quad t >\tau, u_0 \in L^{p}(\Omega_{\varepsilon}), 
	\end{equation*}
	where $M$ and $\tilde{M}$ depend of $N$ and $p$, $\nu$ comes from the exponential decay that the semigroup generated by $-A_{\varepsilon}(\tau)$ possesses. For the case $\varepsilon = 0$,  $t>\tau$ and $u_0 \in L^{p}(\Omega)\times L^{p}(0,1)  $,
	\begin{equation*}
	\norma{U_{0} (t,\tau)  u_0     }_{W^{1,p} (\Omega) \times W^{1,p}(0,1)} \leq M e^{-\nu (t-\tau)} (t-\tau)^{-\frac{1}{2}-(1-\alpha)} \norma{u_0}_{L^{p}(\Omega) \times L^{p}(0,1)},
	%{\small \quad t >0, u_0 \in L^{p}(\Omega)\times L^{p}(0,1) }
	\end{equation*}
	\begin{equation*}
	\norma{U_{0} (t,\tau)  u_0     }_{L^{\infty} (\Omega)\times L^{\infty}(0,1)} \leq \tilde{M} e^{-\nu (t-\tau)} (t-\tau)^{-\frac{1}{2}-(1-\alpha)} \norma{u_0}_{L^{p}(\Omega) \times L^{p}(0,1)}.
	% \quad t >0, u_0 \in L^{p}(\Omega)\times L^{p}(0,1) 
	\end{equation*}
	
	The constants can be chosen uniformly for $\varepsilon \in [0,1]$.
\end{lemma}

We can now prove:

\begin{theorem}\label{t41.1,41.2,41.4}
	For any $\varepsilon \in [0,1]$, let $u_0 \in U_{p}^{\varepsilon}.$ The solution of the semilinear evolution problem 
	$$
	\begin{cases}
	u_t+A_{\varepsilon}(t) u = F_{\varepsilon} (u), \\
	u(\tau) = u_0 \in U_{p}^{\varepsilon},
	\end{cases}
	$$
	is globally defined. Moreover, there exist $K_{\infty}, K_{\frac{1}{2}}>0$, depending only on $N$ and $p$ (they can be chosen uniformly on $\varepsilon$), such that
	\begin{equation}\label{41.I}
	\limsup_{s\rightarrow  - \infty} \norma{u(t,s, u_0)}_{\infty} \leq K_{\infty},			
	\end{equation}
	\begin{equation}\label{41.II}
	\limsup_{s\rightarrow - \infty} \norma{u(t,s,u_0)}_{\frac{1}{2}} \leq K_{\frac{1}{2}}
	\end{equation}
	and those limits are uniform for initial conditions in bounded subsets of $U_{p}^{\varepsilon}$. The norms above are the following
	\[ \norma{\cdot}_{\infty} = \norma{\cdot}_{L^{\infty} (\Omega_{\varepsilon})}   \mbox{ and }   \norma{\cdot}_{\frac{1}{2} }= \norma{\cdot}_{W^{1,p} (\Omega_{\varepsilon})}      ,\mbox{ for }\varepsilon \in (0,1],\]
	\[ \norma{\cdot}_{\infty} = \norma{\cdot}_{L^{\infty} (\Omega)\times L^{\infty}(0,1)}   \mbox{ and }   \norma{\cdot}_{\frac{1}{2} }= \norma{\cdot}_{W^{1,p} (\Omega) \times W^{1,p}(0,1)}      ,\mbox{ for }\varepsilon=0.\]

	The only additional hypothesis we need for the limiting case is that $\alpha +\frac{\delta}{2} >1$.
	%, which is not any restriction, as we have already mentioned.
\begin{proof}
	
	We give a proof for the case $\varepsilon =0$.
	Let $B \subset D(A_0(\cdot))$ be a bounded set in $U_p^{0}$, that is, for $u_0 \in B$, $\norma{u_0}_{U_{p}^{0}} \leq L$. We know from \eqref{40.IV} and from the embedding $D(A_{0}(\cdot)) \hookrightarrow W^{1,p}(\Omega) \oplus W^{1,p}(0,1)$ that $\norma{u(t,s,u_0)}_{\infty} \leq \norma{\gamma(t,s,|u_0|)}_{\infty}.$
	
	From the variation of parameters formula
	$$\gamma (t,s,|u_0|) = U_{0} (t,s)|u_0| + \int_{s}^{t} U_{0}(t,\theta) \left[   \xi \gamma (\theta, s,|u_0|)   +K_M(\theta) \right]					d\theta .$$
	
	By Lemma \ref{l40.3} and Gronwall's Lemma, denoting $\phi(t) = \norma{\gamma(t,s,|u_0|)}_{\infty}$ and $\beta =\frac{1}{2} +(1-\alpha)$, we have 
	%$\beta \in (0,1)$ by hypothesis,
	{\small 		\begin{align*}
		\phi (t)   &  \leq \tilde{M}e^{-\nu (t-s)}(t-s)^{-\beta}\norma{u_0}_{U_p^{\varepsilon}} + \int_{s}^{t} C e^{-\nu (t-\theta)} (t-\theta)^{-\beta} \left[\phi(\theta) + M \right] d\theta \\
		&  \leq \tilde{M}e^{-\nu (t-s)}(t-s)^{-\beta}\norma{u_0}_{U_p^{\varepsilon}} + CM \int_{s}^{t} e^{-\nu (t-\theta)} (t-\theta)^{-\beta} d\theta  +\\
		& \quad + \int_{s}^{t} C e^{-\nu (t-\theta)} (t-\theta)^{-\beta} \left[\phi(\theta) + M \right] d\theta \\
		& \leq \left\{ CLe^{-\nu (t-s)} (t-s)^{-\beta}    + C\int_{s}^{t} e^{-\nu (t-\theta)} (t-\theta)^{-\beta} d\theta \right\} exp \left(C{\int_{s}^{t}  e^{-\nu(t-\theta)} (t-\theta)^{-\beta} d\theta }\right),
		\end{align*}}where $C$ is a constant that can be chosen uniformly for $\varepsilon \in [0,1]$, as a consequence of Lemma \ref{l40.3}. 
	
	The integral appearing after the last inequality can be uniformly bounded in $s$ and $t$, as a matter of fact,
	$$\int_{s}^{t} e^{-\nu (t-\theta)} (t-\theta)^{-\beta} d\theta  \leq \dfrac{e^{-\nu}}{\nu} + \dfrac{1}{1-\beta}= C_{\beta}.$$
	
	%				
	%				Indeed,
	%				\begin{align*}
	%				\int_{s}^{t} e^{-\nu (t-\theta)} (t-\theta)^{-\beta} d\theta &  \leq \int_{-\infty}^{t} e^{-\nu (t-\theta)} (t-\theta)^{-\beta} d\theta  \\
	%				& \leq \int_{-\infty}^{t-1} e^{-\nu (t-\theta)} (t-\theta)^{-\beta} d\theta + \int_{t-1}^{t} e^{-\nu (t-\theta)} (t-\theta)^{-\beta} d\theta \\
	%				& \leq \int_{-\infty}^{t-1} e^{-\nu (t-\theta)}d\theta + \int_{t-1}^{t} (t-\theta)^{-\beta} d\theta \\
	%				& \leq \left[  \dfrac{e^{-\nu (t-\theta)}}{\nu}   \right]_{-\infty}^{t-1} + \left[  \dfrac{(t-\theta)^{1-\beta}}{-(1-\beta ) }   \right]_{t-1}^{t}\\
	%				& \leq  \dfrac{e^{-\nu}}{\nu}    +   \dfrac{1}{(1-\beta ) }
	%				\end{align*}
	%				
	
	Therefore, taking $K_{\infty}-1=CC_{\beta}e^{CC_{\beta}}$, we have
\begin{equation}\label{atracao_uniforme}
	\phi (t) \leq e^{CC_{\beta}} CL e^{-\nu (t-s)} (t-s)^{-\beta}  + K_{\infty}-1.
\end{equation}

	This is clearly bounded for $t$ in bounded intervals (away from the initial time $s$). Besides, noting that $$e^{-\nu (t-s)}   (t-s)^{-\beta}\stackrel{s\rightarrow -\infty}{\longrightarrow} 0,$$
	we have 
	$\limsup_{s\rightarrow  - \infty} \norma{\gamma (t,s,\tau)}_{\infty} =K_{\infty}-1 \leq K_{\infty},$
	uniformly in the bounded set $B$.
	
	\mbox{ }
	
	For the general case where $B$ is a subset of $U_p^{0}$, we note that the nonlinear process regularizes immediately, meaning that if we fix a small $0<\eta<t_0$, 
	$u(s+\eta, s, u_0) \in D(A_0(\cdot)) $
	and, from Proposition \ref{t25.5}, if $B \subset B_{U_p^{0}}[0,L]$, then 
	$$\norma{u(s+\eta,s,u_0)}_{U_p^{0}} \leq \dfrac{CL}{\eta^{1-\alpha}}.$$
	
	Therefore, for $t> s+\eta$, $u(t,s,u_0) = u(t,s+\eta, \tilde{u}_0)$, where $ \tilde{u}_0 = u(s+\eta,s,u_0)$ and $\tilde{u}_0$ belongs to $ B_{U_p^{0}} \left[ 0, \frac{CL}{\eta^{1-\alpha}}  \right] \cap D(A_0(\cdot)).$ Hence
	$$\limsup_{s\rightarrow  - \infty} \norma{u(t,s,u_0)}_{\infty} = \limsup_{s\rightarrow  - \infty} \norma{u(t,s+\eta,\tilde{u}_0)}_{\infty} \leq K_{\infty}. $$
	
	\mbox{ }
	
	We know now that $\norma{u(t,s,u_0)}_{\infty}$ can be uniformly bounded for $u_0 \in B\subset B_{U_p^{0}}[0,L]$, at least for $t$ not close to the initial time $s$. In that case, $F_0(u(\cdot)) \in U_p^{0}$ and we have a bound
	$$\norma{F_0(u(\theta, s, u_0))}_{U_p^{0}} \leq C, $$
	where $C=C( \norma{u(\theta, s, u_0)}_{U_p^{0}}  )$.
	
	Consider $t \geq \theta$, where $\theta = s+a$ and $a$ is large enough so that 
	$$\norma{u(t,s,u_0)}_{\infty} \leq K_{\infty}.$$
	
	We also note that, from the fact that $u(t,s,u_0) \in D(A_0(\cdot)) \hookrightarrow \mathcal{C}(\Omega)\oplus \mathcal{C}(0,1)$, 
	we have $$ \norma{u(t,s,u_0)} _{U_p^{0}} \leq C \norma{u(t,s,u_0)  }_{\infty}  $$ and $C$ can be taken uniformly in $\varepsilon \in [0,1]$.
	
	So, writing $\beta = -\frac{1}{2} - (1-\alpha)$, we have
	$$ \norma{u(t,s,u_0)}_{\frac{1}{2}}  =\norma{u(t,t_0,u(t_0,s,u_0))}_{\frac{1}{2}} \leq Ce^{-\nu (t-t_0)}(t-t_0)^{-\beta}K_{\infty} + C \int_{t_0}^{t}e^{-\nu (t-\theta) }  (t-\theta)^{-\beta} d\theta. $$
	
	The first term goes to zero as $s \rightarrow -\infty$ ($t_0 \rightarrow - \infty$) and the second converges (independent of $t$ or the bounded set $B$). This way we proved the existence of a constant $K_{\frac{1}{2}} $ such that  \eqref{41.II} is satisfied.
\end{proof}
\end{theorem}

%%%%%%%%%%%%%%%%%%%%%%%%%%%%%%%%%%%%%%

%%%%%%%%%%%%%%%%%%%%%%%%%%%%%%%%%%%%%%

\subsection{Existence of pullback attractor}

As pointed out in the previous theorem, the constants $K_{\infty}$ and $K_{\frac{1}{2}}$ are independent of $t$ and $B$, and can even be chosen uniformly in $\varepsilon \in [0,1]$. Thus $$B_{U_{\infty}^{0}}   \left[0,K_{\infty}\right]$$ is a set that pullback absorbs bounded sets of $U_{p}^{0}$ in the $U_{\infty}^{0}$ norm, and $$B_{W^{1,p}(\Omega)\oplus W^{1,p}(0,1)}   \left[0,K_{\frac{1}{2}}\right]$$
is a set that pullback absorbs bounded sets of $U_p^{0}$ in the $W^{1,p}(\Omega)\oplus W^{1,p}(0,1)$ norm.

The existence of those two closed balls that pullback absorb in stronger norms than the one defined for $U_p^{\varepsilon}$ allows us to prove the main result on this section, Theorem \ref{t41.3_41.5}, which states the existence of pullback attractors. Besides, the uniformity of the constants $K_{\infty}$ and $K_{\frac{1}{2}}$ implies that we can obtain uniform bounds for the attractors.

\begin{theorem}\label{t41.3_41.5}
	Let 
	%$\frac{p(2N+1)}{2p+1}<q \leqslant p$ and 
	$N< p$, $X_{\varepsilon} = U_p^{\varepsilon}$
	%, $Y_{\varepsilon}= U_q^{\varepsilon}$ 
	and $f : \R \rightarrow \R$ satisfying both growth and dissipativeness conditions (\eqref{G**} and \eqref{D}). If $K_{\infty}$ and $K_{\frac{1}{2}}$ are the constants obtained previously, then, for every $\varepsilon \in [0,1]$, the semilinear evolution problem
	$$\begin{cases}
	u_t+ A_{\varepsilon}(t) u = F_{\varepsilon}(u),\\
	u(\tau ) = u_0 \in U_p^{\varepsilon},
	\end{cases}$$
	has a pullback attractor $\{\mathcal{A}_{\varepsilon} (t); t\in\R\}$ in $X_{\varepsilon}$ such that 
	\begin{equation}\label{uni_inf}
	\bigcup_{t\in \R} \mathcal{A}_{\varepsilon}(t)  \subset B_{U_{\infty}^{\varepsilon}}\left[  0,K_{\infty}  \right].\end{equation}
	
	Besides, for some $0<\nu <1$,
	\begin{equation} \label{uni_0,5}
	\bigcup_{t\in \R} \mathcal{A}_{\varepsilon}(t)  \subset {C}^{\nu}(\overline{\Omega_{\varepsilon}}), \quad \mbox{ when }\varepsilon \in (0,1],\end{equation}
	\begin{equation*}
	\bigcup_{t\in \R} \mathcal{A}_{0}(t)  \subset {C}^{\nu}(\overline{\Omega}) \oplus {C}^{\nu}(0,1), \quad \mbox{ when }\varepsilon =0.\end{equation*}
	
	Moreover, $\mathcal{A}_{\varepsilon}(t) $  pullback attracts bounded sets of $X_{\varepsilon}$ in the topology of ${C}(\overline{\Omega_{\varepsilon}}) $  and $  {C}(\overline{\Omega}) \oplus {C}([0,1]),$ for $\varepsilon \in(0,1]$ and $\varepsilon =0$, respectively.
	
	\begin{proof}
		Once $W^{1,p}(\Omega_{\varepsilon})$ and $W^{1,p}(\Omega) \oplus W^{1,p} (0,1)$ are compactly embedded in $U_p^{\varepsilon}$ and $U_p^{0}$, respectively, we can conclude that the immersion of the two balls
		$$B_{W^{1,p}(\Omega_{\varepsilon})  } \left[0,K_{\frac{1}{2}}\right]\quad \mbox{ and } \quad B_{W^{1,p}(\Omega)\oplus W^{1,p}(0,1)}   \left[0,K_{\frac{1}{2}}\right]$$ are compact sets of $U_p^{\varepsilon}$ and $U_p^{0}$, respectively, that pullback absorb bounded sets. Therefore, it follows from Theorem \ref{t21.1} that the associated process possesses pulback attractors $\{ \mathcal{A}_{\varepsilon}(t) ; t\in \R \}$, for $\varepsilon \in [0,1]$.
		
		From the fact that $B_{U_{\infty}^{\varepsilon}}\left[ 0,K_{\infty} \right]$ is pullback absorbing, we then have \eqref{uni_inf}. For the assertion \eqref{uni_0,5} of the theorem, we consider the case $\varepsilon \in (0,1]$. The other needs exactly the same arguments.
		
		We know that there exists  a bounded subset in $W^{1,p}(\Omega_{\varepsilon})$ that pullback absorbs. Moreover, this space is continuously embedded in ${C}^{\nu}(  \overline{\Omega_{\varepsilon}}) $ for $0< \nu < 1- \frac{N}{p}$ (see \cite[Theorem 9.16]{Brezis}) and
		$$\bigcup_{t\in \R}  \mathcal{A}_{\varepsilon}(t)   \subset W^{1,p}(\Omega) \hookrightarrow {C}^{\nu}(  \overline{\Omega_{\varepsilon}}).$$
		
		%
		% Aqui eu considerei a realização em ... ao invés de ... pois claramente temos a setorialidade lá e a definição de processo lá
		%
		
		For the last assertion and $\varepsilon =0$, we consider the realization of the operator $A_0(t)$ in ${C}(\overline{\Omega}) \oplus {C}(0,1)$, which is sectorial (Proposition \ref{p22.1}). 
		The nonlinear process generated in this case $\tilde{S}_0(t,s) : {C}(\overline{\Omega}) \oplus {C}(0,1) \rightarrow {C}(\overline{\Omega}) \oplus {C}(0,1)$ will also have a pullback attractor $\tilde{\A}_{0}(t)$. Clearly $\tilde{\A}_{0}(t) \subset \A_0(t)$, since the second one attracts a bigger universe. 
		
		On the other hand, given any bounded set $B_0 \subset U_p^{0}$, $S_0(t,s)B_0$ is bounded in $W^{1,p}(\Omega) \oplus W^{1,p}(0,1)$, for any $t>s$. Hence, $S_0(t,s)B_0$ is bounded in ${C}(\overline{\Omega}) \oplus {C}(0,1)$ and is attracted by $\tilde{\A}_0(t)$ in the topology of ${C}(\overline{\Omega}) \oplus {C}(0,1)$ and, consequently, in the topology of $U_p^{0}$.
		
		Therefore, $\A_0(t) \subset \tilde{\A}_0(t)$, the equality of those sets follows and $\A (t)$ attracts bounded sets of $U_p^{0}$ also in the topology of ${C}(\overline{\Omega}) \oplus C(0,1)$.
	\end{proof}
\end{theorem}

%%%%%%%%%%%%%%%%%%%%%%%%%%%%%%%%%%%%%%

%%%%%%%%%%%%%%%%%%%%%%%%%%%%%%%%%%%%%%

\subsection{Existence of uniform attractor}

	A closer look to the proof of Theorem \ref{t41.1,41.2,41.4} (inequality \eqref{atracao_uniforme}) shows us that $$\norma{u(t,s,u_0)}_{\infty} \leq \norma{\gamma(t,s,|u_0|)}_{\infty} \leq e^{CC_{\beta}} CL e^{-\nu (t-s)} (t-s)^{-\beta}  + K_{\infty}-1$$
and, if $t= s+\tau$, $\tau>0$, 
$$\norma{u(\tau + s,s,u_0)}_{\infty} \leq e^{CC_{\beta}} CL e^{-\nu \tau} (\tau)^{-\beta}  + K_{\infty}-1.$$

This implies that 
$$\limsup_{\tau\rightarrow  \infty} \norma{u(\tau +s ,s,u_0)}_{\infty}  \leq K_{\infty} $$
and this limsup is uniform for $s \in \R $ and $u_0$ in $B_{U_p^{\varepsilon}}[0,L]$. Therefore, the %compact 
pullback absorbing set
$$B_{U_{\infty}^{\varepsilon}}   \left[0,K_{\infty}\right], \quad \varepsilon\in [0,1],$$
is also forward attracting.

 The same analysis can be done to $\norma{u(s+\tau,s,u_0)}_{\frac{1}{2}}$ and we have
$$\limsup_{\tau\rightarrow  \infty} \norma{u(\tau +s ,s,u_0)}_{\frac{1}{2}}  \leq K_{\frac{1}{2}}, $$
uniform for $s\in \R$ and $u_0$ in a bounded set of $U_p^{\varepsilon}$. Hence, the compact pullback absorbing sets
	$$B_{W^{1,p}(\Omega_{\varepsilon})  } \left[0,K_{\frac{1}{2}}\right], \mbox{for }\varepsilon \in (0,1],\mbox{ and }  B_{W^{1,p}(\Omega)\oplus W^{1,p}(0,1)}   \left[0,K_{\frac{1}{2}}\right], \mbox{ for } \varepsilon =0$$
	are also forward attracting. 
	
	This discussion proves the next result.

\begin{proposition}\label{lim_sup_unif}
	Let $K_{\infty}$ and $K_{\frac{1}{2}}$ be the constant obtained in Theorem \ref{t41.1,41.2,41.4}, $B \subset U_p^{\varepsilon}$ a bounded set and $\{S_{\varepsilon} (t,s); t>s \}$ the process generated by the solution of the semilinear evolution equation
		$$
	\begin{cases}
	u_t+A_{\varepsilon}(t) u = F_{\varepsilon} (u), \\
	u(\tau) = u_0 \in U_{p}^{\varepsilon}.
	\end{cases}
	$$ 
	Then, for $\varepsilon \in(0,1],$ 
	
	$$ \lim_{\tau \rightarrow \infty}  \left(  \sup_{s \in \R} dist \left(  S_{\varepsilon} (\tau + s , s )  B , B_{U_{\infty}^{\varepsilon}}  \left[ 0, K_{\infty}   \right] \right)     \right) =0 ,$$
	$$ \lim_{\tau \rightarrow \infty}  \left(  \sup_{s \in \R} dist \left(  S_{\varepsilon} (\tau + s , s )  B , B_{W^{1,p}(\Omega_{\varepsilon})  } \left[0,K_{\frac{1}{2}}\right] \right)     \right) =0 .$$
	
	For the case $\varepsilon = 0$, we have
	$$ \lim_{\tau \rightarrow \infty}  \left(  \sup_{s \in \R} dist \left(  S_{0} (\tau + s , s )  B , B_{U_{\infty}^{0}}  \left[ 0, K_{\infty}   \right] \right)     \right) =0 ,$$
	$$\lim_{\tau \rightarrow \infty}  \left(  \sup_{s \in \R} dist \left(  S_{0} (\tau + s , s )  B ,B_{W^{1,p}(\Omega)\oplus W^{1,p}(0,1)}   \left[0,K_{\frac{1}{2}}\right] \right)     \right) =0 .$$
\end{proposition}

	The above proposition states that the sets $B_{W^{1,p}(\Omega_{\varepsilon})  } \left[0,K_{\frac{1}{2}}\right] $ and $B_{W^{1,p}(\Omega)\oplus W^{1,p}(0,1)}   \left[0,K_{\frac{1}{2}}\right]$ are compact uniformly attracting sets for the process $S_{\varepsilon} (\cdot, \cdot)$ and $S_0(\cdot, \cdot)$, respectively. In this case, the existence of uniform attractor for those process follows as a consequence of Theorem \ref{exist_uni_atr}.
	
	\begin{theorem}\label{unif_attractor}
		Let 
		%$\frac{p(2N+1)}{2p+1}<q \leqslant p$ and 
		$N< p$, $X_{\varepsilon} = U_p^{\varepsilon}$
		%, $Y_{\varepsilon}= U_q^{\varepsilon}$ 
		and $f : \R \rightarrow \R$ satisfying both growth and dissipativeness conditions (\eqref{G**} and \eqref{D}). If $K_{\infty}$ and $K_{\frac{1}{2}}$ are the constants obtained previously, then, for every $\varepsilon \in [0,1]$, the semilinear evolution problem
		$$\begin{cases}
		u_t+ A_{\varepsilon}(t) u = F_{\varepsilon}(u),\\
		u(\tau ) = u_0 \in U_p^{\varepsilon},
		\end{cases}$$
		has a uniform attractor $\overline{\mathcal{A}}_{\varepsilon}$ in $X_{\varepsilon}$ such that 
			$\overline{\mathcal{A}}_{\varepsilon}  \subset B_{U_{\infty}^{\varepsilon}}\left[  0,K_{\infty}  \right].$	Besides, for some $0<\nu <1$, $\overline{\mathcal{A}}_{\varepsilon}  \subset {C}^{\nu}(\overline{\Omega_{\varepsilon}})$, when $\varepsilon \in (0,1]$, 
		$\overline{\mathcal{A}}_{0}  \subset {C}^{\nu}(\overline{\Omega}) \oplus {C}^{\nu}(0,1)$ and, if $\{\A_{\varepsilon}(t);t\in \R\}$ is the pullback attractor, then
		$$ \cup_{t\in \R} \A_{\varepsilon}(t) \subset \overline{\A}_{\varepsilon}, \quad \forall \varepsilon\in [0,1].$$
		
	%	Moreover, $\mathcal{A}_{\varepsilon}(t) $  pullback attracts bounded sets of $X_{\varepsilon}$ in the topology of ${C}(\overline{\Omega_{\varepsilon}}) $  and $  {C}(\overline{\Omega}) \oplus {C}([0,1]),$ respectively for $\varepsilon \in(0,1]$ and $\varepsilon =0$.
		
	\end{theorem}

\begin{remark}
	Note that the pullback attractor and the uniform attractor focus on different features of the problem. The first one deals with the pullback dynamics and is an invariant set. The invariance of the pullback attractor is replaced by the minimality of the uniform attractor, which is a fixed set in the phase space that attracts forward. Even though the pullback and the forward dynamics are, in general, not related, for the reaction-diffusion equations considered in this work, the forward dynamics was obtained from the pullback dynamics. As a matter of fact, we found, via the pullback approach, a set in the phase space $X_{\varepsilon}$ and then showed that this set is also forward attracting. For further discussion about relations between the pullback and uniform attractor for a nonautonomous problem, we recommend \cite{Bortolan_Caraballo_Carvalho_Langa}, \cite{Bortolan_Carvalho_Langa} and \cite{ChepyzhovVishik}.
\end{remark}

%%
%%
%% acm é um comum usado pela matemática.
%%
%%Ver outros padrões de citação em:  https://www.reed.edu/cis/help/LaTeX/bibtexstyles.html
%%
\bibliographystyle{acm}\bibliography{Existence_pullback_attractor_v10}

\end{document}